\documentclass[11pt]{article}
\usepackage{caption}
\usepackage{url}
\usepackage{verbatim}
\usepackage{enumitem}
\usepackage{amsmath}
\allowdisplaybreaks
\usepackage[usenames,dvipsnames]{color}
\usepackage{amsthm}
\usepackage{amsfonts}
\usepackage{graphicx}
\usepackage{nicefrac}
\usepackage{url}
\usepackage{etoolbox}
\usepackage{algorithm}
\usepackage{algpseudocode}
\usepackage{graphicx, amssymb, subcaption,graphics}
\usepackage{tikz}
\usetikzlibrary{positioning}
\usetikzlibrary{calc,3d}
\def\captionfont{\setb@se{11pt}\protect\footnotesize}
\def\captionfont{\protect\footnotesize}
\usepackage{tikz}
\usetikzlibrary{arrows}

\newcommand{\vertiii}[1]{{\left\vert\kern-0.25ex\left\vert\kern-0.25ex\left\vert #1 \right\vert\kern-0.25ex\right\vert\kern-0.25ex\right\vert}}

\def\0{\mbox{\boldmath $0$}}

\newcommand{\nrm}[1]{\left\| #1 \right\|}
\newcommand{\ciptwo}[2]{\left( #1 , #2 \right)}

\newcommand{\eipx}[2]{\left[ #1 , #2 \right]_{\rm x}}
\newcommand{\eipy}[2]{\left[ #1 , #2 \right]_{\rm y}}
\newcommand{\eipz}[2]{\left[ #1 , #2 \right]_{\rm z}}

\newtheorem{thm}{Theorem}[section]
\newtheorem{prop}[thm]{Proposition}
\newtheorem{cor}[thm]{Corollary}

\newtheorem{rmk}[thm]{Remark}
\newtheorem{lem}[thm]{Lemma}

\def\eqref#1{(\ref{#1})}

\newcommand{\nn}{\nonumber}
\newcommand{\hf}{\frac{1}{2}}
\newcommand{\qf}{\frac{1}{4}}
\newcommand{\bx}{\mathbf{x}}
\newcommand{\trh}{\Delta_h}
\newcommand{\nh}{\nabla_h}

 \def\P{\mbox{$\mathsf{P}_h$}}

\begin{document}

\title{Convergence Analysis and Numerical Implementation of a Second Order Numerical Scheme for the Three-Dimensional Phase Field Crystal Equation}

	\author{
Lixiu Dong\thanks{School of Mathematical Sciences, Beijing Normal University, Beijing 100875, P.R. China (201421130036@mail.bnu.edu.cn)}
	\and	
Wenqiang Feng\thanks{Department of Mathematics, The University of Tennessee, Knoxville, TN 37996 (Corresponding Author: wfeng1@utk.edu)}
	\and
Cheng Wang\thanks{Department of Mathematics, The University of Massachusetts, North Dartmouth, MA  02747 (cwang1@umassd.edu)}	
	\and
Steven M. Wise\thanks{Department of Mathematics, The University of Tennessee, Knoxville, TN 37996 (swise@math.utk.edu)} 
        \and
Zhengru Zhang\thanks{School of Mathematical Sciences, Beijing Normal University, Beijing 100875, P.R. China (zrzhang@bnu.edu.cn)} 
}

\maketitle
\numberwithin{equation}{section}

\begin{abstract}
In this paper we analyze and implement a second-order-in-time numerical scheme for the three-dimensional phase field crystal (PFC) equation. The numerical scheme was proposed in \cite{hu09}, with the unique solvability and unconditional energy stability established. However, its convergence analysis remains open. We present a detailed convergence analysis in this article, in which the maximum norm estimate of the numerical solution over grid points plays an essential role. Moreover, we outline the detailed multigrid method to solve the highly nonlinear numerical scheme over a cubic domain, and various three-dimensional numerical results are presented, including the numerical convergence test, complexity test of the multigrid solver and the polycrystal growth simulation. 
\end{abstract}

\section{Introduction}

Defects, such as vacancies, grain boundaries, and dislocations, are observed in crystalline materials, and a precise and accurate understanding of their formation and evolution is of great interest. The phase field crystal (PFC) model was proposed in \cite{elder02} as a new approach to simulate crystal dynamics at the atomic scale in space but on diffusive scales in time. This model naturally incorporates elastic and plastic deformations, multiple crystal orientations and defects and has already been used to simulate a wide variety of microstructures, such as epitaxial thin film growth~\cite{elder04}, grain growth~\cite{stefanovic06}, eutectic solidification~\cite{elder07}, and dislocation formation and motion~\cite{provatas07, stefanovic06}.  
The idea is that the phase variable describes a coarse-grained temporal average of the number density of atoms, and the approach can be related to dynamic density functional theory \cite{backofen07, marconi99}.  The method represents a significant advantage over other atomistic methods, such as molecular dynamics methods where the time steps are constrained by atomic-vibration time scales.	 More detailed descriptions are available in~\cite{Achim06,Achim09,aland2012particles,backofen2014phase, Berry08, Berry11, elder10a, Heinonen14, Mellenthin08,praetorius2011phase,praetorius2015navier, Ramos08, Ramos10, stefanovic06, stefanovic09, wheeler06, wu10, wu09,yu2011morphological}, and the related works for the amplitude expansion approach could be found in~\cite{Athreya06, elder10a, Goldenfeld05, Goldenfeld06, guan16a, Spatschek10, Yeon10}. 

Consider the dimensionless energy of the form~\cite{elder04, elder02, swift77}: 
\begin{equation}
 E(\phi)=\int_\Omega\left[\qf\phi^4+\frac{1-\epsilon}{2}\phi^2-|\nabla\phi|^2 +\hf(\Delta\phi)^2\right]d\bx, \label{energy-PFC} 
\end{equation}
where $\Omega = (0,L_x)\times(0,L_y) \times(0,L_z)\subset \mathbb{R}^3$, $\phi:\Omega\rightarrow \mathbb{R}$ is the atom density field, and $\epsilon > 0$ is a constant. We assume that $\phi$ is periodic on $\Omega$. This model naturally incorporates elastic and plastic deformation of the crystal and the various crystal defects. The PFC equation~\cite{elder04, elder02} is given by the $H^{-1}$ gradient flow associated with the energy (\ref{energy-PFC}): 
\begin{equation} \label{equation-PFC}
\begin{aligned}
&\phi_t=\nabla\cdot\left(M(\phi)\nabla\mu\right) ,\quad \text{in}\quad \Omega_T= \Omega\times(0,T) ,\\
&\mu= \delta_\phi E = \phi^3+(1-\epsilon)\phi+2\Delta \phi+\Delta^2 \phi ,\quad \text{in}\quad \Omega_{T},\\
&\phi(x,y,z,0)= \phi_0 (x,y,z),\quad \text{in}\quad \Omega. 
\end{aligned}
\end{equation}
in which $M(\phi)>0$ is a mobility, $\mu$ is the chemical potential. Periodic boundary conditions are imposed for $\phi$, $\Delta \phi$ and $\mu$. 

The PFC equation is a high-order (sixth-order) nonlinear partial differential equation. There have been some related works to develop numerical schemes for the PFC equation.  Cheng and Warren~\cite{cheng08} introduced a linearized spectral scheme, similar to one for the Cahn-Hilliard equation analyzed in~\cite{vollmayerlee03}. This scheme is not expected to be provably unconditionally energy stable.  The finite element PFC method of Backofen~\emph{et al.}~\cite{backofen07} employs what is essentially a standard backward Euler scheme, but where the nonlinear term $\phi^3$ in the chemical potential is linearized via $(\phi^{k+1})^3 \approx 3(\phi^k)^2\phi^{k+1} - 2(\phi^k)^3$.  Both energy stability and solvability are issues for this scheme, because the term $2\Delta\phi$ is implicit in the chemical potential. Tegze~\emph{et al.}~\cite{tegze09} developed a semi-implicit spectral scheme for the binary PFC equations that is not expected to unconditionally stable. Also see other related numerical works~\cite{Athreya07, cao2015two, guoR16, Hirouchi09} in recent years. 

The energy stability of a numerical scheme has always been a very important issue, since it plays an essential role in the accuracy of long time numerical simulation.  The standard convex splitting scheme, originated from Eyre's work~\cite{eyre98}, has been a well-known approach to achieve numerical energy stability. In this framework, the convex part of the chemical potential is treated  implicitly, while the concave part is updated explicitly. A careful analysis leads to the unique solvability and unconditional energy stability of the numerical scheme, unconditionally with respect to the time and space step sizes. Such an idea has been applied to a wide class of gradient flows in recent years, and both first and second order accurate in time algorithms have been developed. See the related works for the PFC equation and the modified PFC (MPFC) equation~\cite{baskaran13a, baskaran13b,bueno2016three,grasselli2015energy, hu09, wang10c, wang11a, wise09a};  
the epitaxial thin film growth models~\cite{chen12, chen14, shen12, wang10a}; Cahn-Hilliard equation~\cite{diegel16, guo16}; non-local Cahn-Hilliard-type models~\cite{guan14b, guan14a}, the Cahn-Hilliard-Hele-Shaw (CHHS) and related models \cite{chen2016chhs,chen16, collins13, diegel15a, feng12, wise10}, etc. 

On the other hand, a well-known drawback of the first order convex splitting approach is that an extra dissipation has been added to ensure unconditional stability; in turn, the first order numerical approach introduces a significant amount of numerical error~\cite{christlieb14}. For this reason, second-order energy stable methods have been highly desirable. 

There have been other related works of ``energy stable" schemes for the certain gradient flows in recent years. For example, an alternate variable is used in~\cite{guillen14}, denoted as a second order approximation to $v=\phi^2-1$ in the Cahn-Hilliard model. A linearized, second order accurate scheme is derived as the outcome of this idea, and the stability is established for a modified energy. A similar idea has been applied to the PFC model in a more recent article~\cite{yang16b}. However, such an energy stability is applied to a pair of numerical variables $(\phi, v)$, and an $H^2$ stability for the original physical variable $\phi$ has not been justified. As a result, the convergence analysis is not available for this numerical approach. Similar methodology has been reported in the invariant energy quadratization (IEQ) approach~\cite{han16a, yang16a, zhao16a}, etc.  

In comparison, a second order numerical scheme was proposed and studied for the PFC equation in~\cite{hu09}. By a careful choice of the second order temporal approximations to each term in the chemical potential, the unique solvability and unconditional energy were justified at a theoretical level, with the centered difference discretization taken in space. In particular, this energy stability is derived with respect to the original energy functional, combined with an auxiliary, non-negative correction term, so that a uniform in time $H^2$ bound is available for the numerical solution. Meanwhile, a detailed convergence analysis has not been theoretically reported for the proposed second order scheme, although the full convergence order was extensively demonstrated in the numerical experiments. The key difficulty in the convergence analysis is associated with the maximum norm bound estimate for the numerical solution, and such a bound plays an essential role in the theoretical convergence derivation. In more details, the unconditional energy stability indicates a uniform in time $H^2$ bound of the numerical solution at a discrete level. Although the Sobolev embedding from $H^2$ to $L^\infty$ is straightforward in three-dimensional space, a direct estimate for the corresponding grid function is not directly available. In two-dimensional space, such a discrete Sobolev embedding has been proved in the earlier works~\cite{hu09, wise09a}, using a complicated calculations of the difference operators. However, as stated in Remark 12 of \cite{hu09}, {\em ``the proof presented in~\cite{wise09a} does not automatically extend to three dimensions. This is because a discrete Sobolev inequality is used to translate energy stability into point-wise stability, and the inequality fails in three dimensions. We are currently studying the three dimensional case in further detail."} 

In this paper, we provide a detailed convergence analysis for the fully discrete scheme formulated in \cite{hu09}, which is shown to be second order accurate in both time and space. In particular, the maximum norm estimate of the three-dimensional numerical solution is accomplished via a discrete Fourier transformation over a uniform numerical grid, so that the discrete Parseval equality is valid. And also, the equivalence between the discrete and continuous $H^2$ norms for the numerical grid function and its continuous version, respectively, can be established. In turn, the discrete Sobolev inequality is obtained from its continuous version. Such an $\ell^\infty$ bound of the discrete numerical solution is crucial, so that the convergence analysis could go through for the scheme. Moreover, the Crank-Nicholson approximation to the surface diffusion term poses another challenge in the convergence proof, since the diffusion coefficients at time steps $t^{k+1}$ and $t^k$ are equally distributed, in comparison with an alternate second order approximation reported in a few recent works~\cite{diegel16, guo16}, in which the diffusion coefficients are distributed at time steps $t^{k+1}$ and $t^{k-1}$, respectively. To overcome this difficulty, we have to perform an error analysis at time instant $t^{k+1/2}$, in combination a subtle estimate for the numerical error in the concave diffusion term. 

In addition, we also present various numerical simulation results of three-dimensional PFC model in this article. It is noted that most numerical results for the PFC equation reported in the existing literature are two-dimensional, or over a two-dimensional surface; see~\cite{baskaran13a, cheng08, dehghan16, gomez2012unconditionally, hu09, zhang13}, etc. For a gradient flow in which the nonlinear terms takes a form of $\phi^3$ pattern, great efficiency and accuracy of the nonlinear multi-grid solver have been extensively demonstrated in the numerical experiments; see the related works~\cite{baskaran13a, collins13, guo16, hu09, guan14b, guan14a, wise10}. We apply the nonlinear multi-grid algorithm to implement the three-dimensional numerical scheme; its great numerical efficiency enables us to compute the three-dimensional model using local servers. Both the numerical accuracy check and the detailed numerical simulation results of three-dimensional  polycrystal growth are reported. 

The rest of paper is organized as follows. In Section~\ref{sec:space}, we introduce the finite difference spatial discretization in three-dimensional space, and review a few preliminary estimates.  In Section~\ref{sec:numerical scheme} we review the second order numerical scheme proposed in~\cite{hu09}, and state the main theoretical results. The detailed convergence proof is given by Section~\ref{sec:convergence}. Furthermore, the details of the three-dimensional multi-grid solver is outlined in Section~\ref{sec:multigrid}.  Subsequently, the numerical results are presented in Section~\ref{sec:numerical results}. Finally, some concluding remarks are made in Section~\ref{sec:conclusion}.

\section{Finite difference discretization and a few preliminary estimates} \label{sec:space} 

For simplicity of presentation, we denote $( \cdot , \cdot )$ as the standard $L^2$ inner product, and $\| \cdot \|$ as the standard $L^2$ norm, and $\| \cdot \|_{H^m}$ as the standard $H^m$ norm. We use the notation and results for some discrete functions and 
operators from~\cite{guo16,wise10, wise09}. 
Let $\Omega = (0,L_x)\times(0,L_y)\times(0,L_z)$, where for simplicity, we assume $L_x =L_y=L_z =: L > 0$. It is also assumed that $h_x = h_y = h_y =h$ and we denote $L = m\cdot h$, where $m$ is a positive integer. The parameter $h = \frac{L}{m}$ is called the mesh or grid spacing. We define the following two uniform, infinite grids with grid spacing $h>0$:
	\[
E := \{ x_{i+\hf} \ |\ i\in {\mathbb{Z}}\}, \quad C := \{ x_i \ |\ i\in {\mathbb{Z}}\},
	\]
where $x_i = x(i) := (i-\hf)\cdot h$. Consider the following 3D discrete periodic function spaces: 
	\begin{eqnarray*}
	\begin{aligned}
{\mathcal C}_{\rm per} &:= \left\{\nu: C\times C
\times C\rightarrow {\mathbb{R}}\ \middle| \ \nu_{i,j,k} = \nu_{i+\alpha m,j+\beta m, k+\gamma m}, \ \forall \, i,j,k,\alpha,\beta,\gamma\in \mathbb{Z} \right\},
	\\
{\mathcal E}^{\rm x}_{\rm per} &:=\left\{\nu: E\times C\times C\rightarrow {\mathbb{R}}\ \middle| \ \nu_{i+\frac12,j,k}= \nu_{i+\frac12+\alpha m,j+\beta m, k+\gamma m}, \ \forall \, i,j,k,\alpha,\beta,\gamma\in \mathbb{Z}\right\} .
	\end{aligned}
	\end{eqnarray*}		
The spaces  ${\mathcal E}^{\rm y}_{\rm per}$ and ${\mathcal E}^{\rm z}_{\rm per}$ are analogously defined. The functions of ${\mathcal C}_{\rm per}$ are called {\emph{cell centered functions}}. The functions of ${\mathcal E}^{\rm x}_{\rm per}$, ${\mathcal E}^{\rm y}_{\rm per}$, and ${\mathcal E}^{\rm z}_{\rm per}$,  are called {\emph{east-west face-centered functions}},  {\emph{north-south face-centered functions}}, and {\emph{up-down face-centered functions}}, respectively.  We also define the mean zero space 
	\[
\mathring{\mathcal C}_{\rm per}:=\left\{\nu\in {\mathcal C}_{\rm per} \ \middle| \overline{\nu} :=  \frac{h^3}{| \Omega|} \sum_{i,j,k=1}^m \nu_{i,j,k}  = 0\right\} .
	\]

We now introduce the important difference and average operators on the spaces:  
	\begin{eqnarray*}
&& A_x \nu_{i+\hf,j,k} := \frac{1}{2}\left(\nu_{i+1,j,k} + \nu_{i,j,k} \right), \quad D_x \nu_{i+\hf,j,k} := \frac{1}{h}\left(\nu_{i+1,j,k} - \nu_{i,j,k} \right),
	\\
&& A_y \nu_{i,j+\hf,k} := \frac{1}{2}\left(\nu_{i,j+1,k} + \nu_{i,j,k} \right), \quad D_y \nu_{i,j+\hf,k} := \frac{1}{h}\left(\nu_{i,j+1,k} - \nu_{i,j,k} \right) , 
\\
&& A_z \nu_{i,j,k+\hf} := \frac{1}{2}\left(\nu_{i,j,k+1} + \nu_{i,j,k} \right), \quad D_z \nu_{i,j,k+\hf} := \frac{1}{h}\left(\nu_{i,j,k+1} - \nu_{i,j,k} \right) , 
	\end{eqnarray*}
with $A_x,\, D_x: {\mathcal C}_{\rm per}\rightarrow{\mathcal E}_{\rm per}^{\rm x}$, $A_y,\, D_y: {\mathcal C}_{\rm per}\rightarrow{\mathcal E}_{\rm per}^{\rm y}$, $A_z,\, D_z: {\mathcal C}_{\rm per}\rightarrow{\mathcal E}_{\rm per}^{\rm z}$. 
Likewise,
	\begin{eqnarray*}
&& a_x \nu_{i, j, k} := \frac{1}{2}\left(\nu_{i+\hf, j, k} + \nu_{i-\hf, j, k} \right),	 \quad d_x \nu_{i, j, k} := \frac{1}{h}\left(\nu_{i+\hf, j, k} - \nu_{i-\hf, j, k} \right),
	\\
&& a_y \nu_{i,j, k} := \frac{1}{2}\left(\nu_{i,j+\hf, k} + \nu_{i,j-\hf, k} \right),	 \quad d_y \nu_{i,j, k} := \frac{1}{h}\left(\nu_{i,j+\hf, k} - \nu_{i,j-\hf, k} \right),
	\\
&& a_z \nu_{i,j,k} := \frac{1}{2}\left(\nu_{i, j,k+\hf} + \nu_{i, j, k-\hf} \right),	 \quad d_z \nu_{i,j, k} := \frac{1}{h}\left(\nu_{i, j,k+\hf} - \nu_{i, j,k-\hf} \right),
	\end{eqnarray*}
with $a_x,\, d_x : {\mathcal E}_{\rm per}^{\rm x}\rightarrow{\mathcal C}_{\rm per}$, $a_y,\, d_y : {\mathcal E}_{\rm per}^{\rm y}\rightarrow{\mathcal C}_{\rm per}$, and $a_z,\, d_z : {\mathcal E}_{\rm per}^{\rm z}\rightarrow{\mathcal C}_{\rm per}$. 
The standard 3D discrete Laplacian, $\Delta_h : {\mathcal C}_{\rm per}\rightarrow{\mathcal C}_{\rm per}$, is given by 
	\begin{equation*}
	\begin{aligned}
\Delta_h \nu_{i,j,k} &:= d_x(D_x \nu)_{i,j,k} + d_y(D_y \nu)_{i,j,k}+d_z(D_z \nu)_{i,j,k}\\
 &= \frac{1}{h^2}\left( \nu_{i+1,j,k}+\nu_{i-1,j,k}+\nu_{i,j+1,k}+\nu_{i,j-1,k}+\nu_{i,j,k+1}+\nu_{i,j,k-1} - 6\nu_{i,j,k}\right).
 \end{aligned}
	\end{equation*}

Now we are ready to define the following grid inner products:  
	\begin{equation*}
	\begin{aligned}
\ciptwo{\nu}{\xi}_2 &:= h^3\sum_{i,j,k=1}^m  \nu_{i,j,k}\xi_{i,j,k},\quad \nu,\, \xi\in {\mathcal C}_{\rm per},\quad
& \eipx{\nu}{\xi} := \ciptwo{a_x(\nu\xi)}{1}_2 ,\quad \nu,\, \xi\in{\mathcal E}^{\rm x}_{\rm per},
\\
\eipy{\nu}{\xi} &:= \ciptwo{a_y(\nu\xi)}{1}_2 ,\quad \nu,\, \xi\in{\mathcal E}^{\rm y}_{\rm per},\quad
&\eipz{\nu}{\xi} := \ciptwo{a_z(\nu\xi)}{1}_2 ,\quad \nu,\, \xi\in{\mathcal E}^{\rm z}_{\rm per}.
	\end{aligned}
	\end{equation*}	

We now define the following norms for cell-centered functions. If $\nu\in {\mathcal C}_{\rm per}$, then $\nrm{\nu}_2^2 := \ciptwo{\nu}{\nu}_2$; $\nrm{\nu}_p^p := \ciptwo{|\nu|^p}{1}_2$ ($1\le p< \infty$), and $\nrm{\nu}_\infty := \max_{1\le i,j,k\le m}\left|\nu_{i,j,k}\right|$.
Similarly, we define the gradient norms: for $\nu\in{\mathcal C}_{\rm per}$,
	\[
\nrm{ \nabla_h \nu}_2^2 : = \eipx{D_x\nu}{D_x\nu} + \eipy{D_y\nu}{D_y\nu} +\eipz{D_z\nu}{D_z\nu}.
	\]
Consequently,
	\[
\nrm{\nu}_{2,2}^2 : =  \nrm{\nu}_2^2+ \nrm{ \nabla_h \nu}_2^2  + \nrm{ \Delta_h \nu}_2^2.
	\]

In addition, the discrete energy $   F_h (\phi): {\mathcal C}_{\rm per}\to \mathbb{R}$ is defined as 
\begin{equation} 
  F_h (\phi) = \frac14 \| \phi \|_4^4 + \frac{1 - \epsilon}{2} \| \phi \|_2^2 - \| \nabla_h \phi \|_2^2 
  + \frac12 \| \Delta_h \phi \|_2^2 . \label{PFC energy-discrete-1} 
\end{equation} 

The following preliminary estimates are cited from earlier works. For more details we refer the reader to \cite{hu09, wise09a}.

\begin{lem}  \label{lemma1}
For any $f, g \in {\mathcal C}_{\rm per}$, the following summation by parts formulas are valid: 
\begin{eqnarray} 
  ( f , \Delta_h g ) = - ( \nabla_h f , \nabla_h g ) ,  \quad 
  ( f , \Delta_h^2 g ) =  ( \Delta_h f , \Delta_h g ) , \quad 
  ( f , \Delta_h^3 g ) = - ( \nabla_h \Delta_h f , \nabla_h \Delta_h g )  . 
  \label{lemma 1-0} 
\end{eqnarray} 
\end{lem} 

\begin{lem}
\label{lemma2}
Suppose $\phi \in {\mathcal C}_{\rm per}$. 
Then
\begin{equation} 
\|\trh \phi\|_2^2 \leq \frac{1}{3\alpha^2}\|\phi\|_2^2+\frac{2\alpha}{3}\|\nh(\trh \phi)\|_2^2 ,  
\label{lemma 2-0} 
\end{equation} 
is valid for arbitrary $\alpha >0$. 
\end{lem}

\begin{lem}
\label{lemma3}
For $\phi \in {\mathcal C}_{\rm per}$, we have the estimate 
\begin{equation} 
 F_h(\phi)\geq  C \|\phi\|_{2,2}^2-\frac{L^3}{4} , \label{lemma 3-0} 
\end{equation} 
with $C$ only dependent on $\Omega$, and $F_h (\phi)$ given by (\ref{PFC energy-discrete-1}). 
\end{lem}

\section{The fully discrete second order numerical scheme and the main results} \label{sec:numerical scheme}
  
Let $N_t \in\mathbb{Z}^+$, and set $\tau:=T/N_t$, where $T$ is the final time.  For our present and future use, we define the canonical grid projection operator $\P: C^0(\Omega)\to {\mathcal C}_{\rm per}$ via $[\P v]_{i,j,k}=v(\xi_i,\xi_j,\xi_k)$. Set $u_{h,s}:=\P u(\cdot, s)$. Then $F_h(u_{h,0},u_{h,s}) + \frac12 \| \nabla_h ( u_{h,s} - u_{h,0} ) \|_2^2 \to F_h(u(\cdot, 0))$ as $h\to 0$ and $s\to 0$ for sufficiently regular $u$. We denote $\phi_{e}$ as the exact solution to the PFC equation \eqref{equation-PFC} and take $\Phi_{i,j,k}^\ell = \P\phi_e (\cdot, t_\ell)$.

Our second order numerical scheme in \cite{hu09} can be formulated as follows: for $1\le \kappa\le N_t$, given $\phi^\kappa,\phi^{\kappa-1}\in {\mathcal C}_{\rm per}$, find $\phi^{\kappa+1},\mu^{\kappa+\hf},\omega^{\kappa+\hf}\in {\mathcal C}_{\rm per}$ periodic such that  
\begin{equation}
\begin{aligned}
&\frac{\phi^{\kappa+1}-\phi^\kappa}{\tau}=\trh\mu^{\kappa+\hf},\\
&\mu^{\kappa+\hf}:=(\phi^{\kappa+\hf})(\phi^2)^{\kappa+\hf}+(1-\epsilon)\phi^{\kappa+\hf}+3\trh\phi^\kappa-\trh\phi^{\kappa-1}+\trh\omega^{\kappa+\hf},\label{1.10}\\
& \omega^{\kappa+\hf}:=\trh\phi^{\kappa+1},
\end{aligned}
\end{equation}
where $\phi^{0}:=\Phi^0$, $\phi^1 := \Phi^1$.

The unique solvability and energy stability have already been established in \cite{hu09}; see the following result. 

\begin{prop}\label{prop:hu09} 
\cite{hu09} 
Suppose that the initial profiles $\phi^0, \phi^1 \in {\mathcal C}_{\rm per}$ satisfy periodic boundary condition, with sufficient regularity assumption for the exact solution $\phi_e$. Given any $(\phi^{m-1}, \phi^m) \in {\mathcal C}_{\rm per}$, there is a unique solution $\phi^{m+1} \in {\mathcal C}_{\rm per}$ to the scheme (\ref{1.10}). The scheme \eqref{1.10}, with starting values $\phi^{0}$ and $\phi^1$, is unconditionally energy stable, \emph{i.e.}, for any $\tau > 0$ and $h>0$, and any positive integer $1\le \kappa\le N_t$,
\begin{equation}\label{ieq:energy dissipative}
F_h (\phi^\kappa) \le F_h (\phi^1) + \frac12 \| \nabla_h ( \phi^1 - \phi^0 ) \|_2^2 \leq   C_0 , 
\end{equation}
in which $C_0$ is independent on $h$, $\tau$, $\epsilon$ and $T$. 
\end{prop}

The $\| \cdot \|_\infty$ bound of a grid function could be controlled with the help of a discrete Sobolev inequality, as stated by the following theorem; its proof will be given in Section~\ref{sec:convergence}. 

\begin{thm} \label{thm:max norm 1} 
 Let $\phi \in {\mathcal C}_{\rm per}$. Then there exists a constant $C$ independent of $\tau$ or $h$ such that
\begin{equation} 
  \|\phi\|_\infty  \leq C \|\phi\|_{2,2} . \label{thm 2-0}   
\end{equation} 
\end{thm}

As a combination of Proposition~\ref{prop:hu09}, Theorem~\ref{thm:max norm 1} and inequality (\ref{lemma 3-0}) in Lemma~\ref{lemma3} , the following $\| \cdot \|_\infty$ estimate for the numerical solution is available. 

\begin{cor} \label{cor:max norm 2}  
For the numerical scheme (\ref{1.10}), we have 
\begin{equation}
   \| \phi^\kappa \|_\infty \le C \left( C_0 + \frac{L^3}{4} \right) := \tilde{C}_0 , \quad \forall \kappa \ge 0 .  \label{cor 3-0}    
\end{equation}
\end{cor}

The main theoretical result is stated in the following theorem; Its proof will be given in Section~\ref{sec:convergence}.

\begin{thm} \label{thm:convergence} 
Suppose the unique solution $\phi_e$ for the three-dimensional PFC equation (\ref{equation-PFC}), with $M (\phi) \equiv 1$, is of regularity class 
                \begin{equation}
			\phi_e \in \mathcal{R} := H^3 (0,T; C_{\rm per}^0) \cap H^2 (0,T; C_{\rm per}^4) \cap L^\infty (0,T; C_{\rm per}^8),
			\label{assumption:regularity.1}
		\end{equation}	
and the initial data $\phi^{0},\phi^{1} \in{\mathcal C}_{\rm per}$ are defined as above. Define $e_{ijk}^\kappa:=\Phi_{ijk}^\kappa-\phi_{ijk}^\kappa$. 		
    	Then, provided $\tau$ and $h$ are sufficiently small, 
    	for all positive integers $\kappa$, such that $\tau\cdot \kappa \le T$, we have
    	\begin{equation}
    	  \|e^{\kappa}\|_2 \leq C(h^2+\tau^2), 
	  \label{convergence-0}
\end{equation}
for some $C>0$ that is independent of $h$ and $\tau$. 
\end{thm}

\section{The detailed convergence analysis} \label{sec:convergence} 

\subsection{The proof of Theorem~\ref{thm:max norm 1}}
We begin with the proof of Theorem~\ref{thm:max norm 1}, which provides a tool to bound the $\| \cdot \|_\infty$ norm of a grid function in terms of its discrete $\| \cdot \|_{2,2}$ norm. 

\begin{proof} 
For a function $\phi\in\mathcal{C}_{\rm per}$  with value $\phi_{ijk}$ at $(x_i,y_j,z_k),$ the IDFT is given by \cite{trefethen2000spectral}: 
\begin{equation}
\phi_{ijk}= \sum_{r,s,t=-R}^{R} \hat{\phi}_{rst} {\rm e}^{2 \pi i (r x_i+ s y_j +t z_k)/L}, \quad i,j,k = 1, \cdots , m , 
\end{equation}
in which $m=2R+1$. In turn, the corresponding interpolation function is defined as 
\begin{equation}
\begin{aligned}
\phi_F(x,y,z)= \sum_{r,s,t=-R}^{R} \hat{\phi}_{rst} {\rm e}^{2 \pi i (r x+ s y +t z)/L} .\nn
\end{aligned}
\end{equation}
Using the Parseval's identity (at both the discrete and continuous levels), we have
$$
  \|\phi\|_2^2=\|\phi_F\|_{L^2}^2=L^3  \sum_{r,s,t=-R}^{R} \left| \hat{\phi}_{rst} \right|^2 , 
$$

\begin{equation}
\begin{aligned}
D_x \phi_{i+\hf,j,k}&=\frac{1}{h}(\phi_{i+1,j,k}-\phi_{i,j,k})\\
                   &=\frac{1}{h} \sum_{r,s,t=-R}^{R} \hat{\phi}_{rst}  \left[ {\rm e}^{2 \pi i ( r x_{i+1} + s y_j+ t z_k )/L} - {\rm e}^{2 \pi i ( r x_i + s y_j + t z_k )/L}\right]\\
                   &=\frac{1}{h} \sum_{r,s,t=-R}^{R} \hat{\phi}_{rst} {\rm e}^{2 \pi i ( r x_{i+\hf} + s y_j + t z_k )/L} \cdot 2i\sin \frac{\pi r h}{L}\\
                   &= \sum_{r,s,t=-R}^{R} u_r \hat{\phi}_{rst} {\rm e}^{2 \pi i ( r x_{i+\hf} + s y_j + t z_k )/L} ,\nn
\end{aligned}
\end{equation}

\begin{equation}
\partial_x \phi_F(x,y,z)= \sum_{r,s,t=-R}^{R} v_r \hat{\phi}_{rst} {\rm e}^{2 \pi i ( r x + s y + t z )/L} , \nn
\end{equation}
with 
\begin{equation}
\begin{aligned}
&u_r=\frac{2i\sin \frac{\pi r h}{L}}{h}, &
   &v_r=\frac{2i \pi r}{L}. \nn
\end{aligned}
\end{equation}
A comparison of Fourier eigenvalues between $|u_r|$ and $|v_r|$ shows that
\begin{equation}
\begin{aligned}
&\frac{2}{\pi}|v_r| \leq |u_r| \leq |v_r|, &
   &-R \leq r\leq R , \nn
\end{aligned}
\end{equation}
$$
[D_x\phi,D_x\phi]_{x}=h^3 \sum_{i,j,k=1}^m |D_x\phi_{i+\hf,j,k}|^2=L^3  \sum_{r,s,t=-R}^{R} |u_r|^2 | \hat{\phi}_{rst} |^2,
$$
$$
\|\partial_x \phi_F\|^2= L^3  \sum_{r,s,t=-R}^{R} |v_r|^2 | \hat{\phi}_{rst} |^2 . 
$$
Then we get 
$$\|\partial_x \phi_F\| ^2=|\frac{v_r}{u_r}|^2[D_x\phi,D_x\phi]_{x} \leq (\frac{\pi}{2})^2[D_x\phi,D_x\phi]_{x},
$$
Similarly, the following estimates are available: 
$$
\|\partial_y \phi_F\| ^2 \leq (\frac{\pi}{2})^2[D_y\phi,D_y\phi]_{y} , 
$$
$$
\|\partial_z \phi_F\| ^2 \leq (\frac{\pi}{2})^2[D_z\phi,D_z\phi]_{z} .
$$
For the second order derivatives, the following estimates are valid: 
\begin{equation}
\begin{aligned}
&d_x(D_x\phi)_{i,j,k}\\
&=\frac{1}{h}(D_x\phi_{i+\hf,j,k}-D_x\phi_{i-\hf,j,k})\\
&=\frac{1}{h}\left( \sum_{r,s,t=-R}^{R}  u_r \hat{\phi}_{rst} {\rm e}^{2 \pi i ( r x_{i+\hf} + s y_j + t z_k )/L} - \sum_{r,s,t=-R}^{R}  u_r \hat{\phi}_{rst} {\rm e}^{2 \pi i ( r x_{i-\hf} + s y_j + t z_k )/L} 
\right)\\
&= \sum_{r,s,t=-R}^{R}  u_r^2 \hat{\phi}_{rst} {\rm e}^{2 \pi i (r x_i + s y_j + t z_k )/L} , \nn
\end{aligned}
\end{equation}

\begin{equation}
\partial_x^2 \phi_F(x,y,z)= \sum_{r,s,t=-R}^{R}  v_r^2 \hat{\phi}_{rst} {\rm e}^{2 \pi i (r x + s y + t z )/L} .\nn
\end{equation}
As a consequence, these inequalities yield the following result: 
$$
\|\partial_x^2 \phi_F\| ^2=|\frac{v_r}{u_r}|^4 \|d_x(D_x\phi)\|_2^2 \leq (\frac{\pi}{2})^4\|d_x(D_x\phi)\|_2^2. 
$$
Similarly, 
$$\|\partial_y^2 \phi_F\|  ^2 \leq (\frac{\pi}{2})^4\|d_y(D_y\phi)\|_2^2,\quad \|\partial_z^2 \phi_F\|  ^2 \leq (\frac{\pi}{2})^4\|d_z(D_z\phi)\|_2^2,$$
and
$$\|\partial_x\partial_y \phi_F\|  ^2 \leq \hf(\|\partial_x^2 \phi_F\|  ^2+\|\partial_y^2 \phi_F\|  ^2),$$
$$\|\partial_x\partial_z \phi_F\|  ^2 \leq \hf(\|\partial_x^2 \phi_F\|  ^2+\|\partial_z^2 \phi_F\|  ^2),$$
$$\|\partial_y\partial_z \phi_F\|  ^2 \leq \hf(\|\partial_y^2 \phi_F\|  ^2+\|\partial_z^2 \phi_F\|  ^2).$$

Then we arrive at 
\begin{equation}
\begin{aligned}
\|\phi_F\|_{H^2}^2&=|\phi_F|_{0,2}^2+|\phi_F|_{1,2}^2+|\phi_F|_{2,2}^2\\
                  &\leq\|\phi_F\|  ^2+\|\partial_x \phi_F\|  ^2+\|\partial_y \phi_F\|  ^2+\|\partial_z \phi_F\|  ^2+2(\|\partial_x^2 \phi_F\|  ^2+\|\partial_y^2 \phi_F\|  ^2+\|\partial_z^2 \phi_F\|  ^2)\\
                  &\leq \|\phi\|_2^2+(\frac{\pi}{2})^2\left([D_x\phi,D_x\phi]_{x}+[D_y\phi,D_y\phi]_{y}+[D_z\phi,D_z\phi]_{z}\right)\\
                  &+2\times(\frac{\pi}{2})^4\left(\|d_x(D_x\phi)\|_2^2+\|d_y(D_y\phi)\|_2^2+\|d_z(D_z\phi)\|_2^2\right)\\
                  &=\|\phi\|_2^2+(\frac{\pi}{2})^2\|\nh \phi\|_2^2+2\times(\frac{\pi}{2})^4\|\trh \phi\|_2^2\\
                  &\leq2 \cdot (\frac{\pi}{2})^4(\|\phi\|_2^2+\|\nh \phi\|_2^2+\|\trh \phi\|_2^2)\\
                  &=2 \cdot (\frac{\pi}{2})^4\|\phi\|_{2,2}^2. \nn
\end{aligned}
\end{equation}
Meanwhile, since $\phi$ is the discrete interpolate of the continuous function, an obvious observation that $\| \phi \|_\infty \le \| \phi_F \|_{L^\infty}$ implies the following estimate: 
$$
 \|\phi\|_\infty \leq \|\phi_F\|_{L^\infty} \leq C_1 \|\phi_F\|_{H^2} \leq \frac{\sqrt{2}\pi^2}{4} C_1 \|\phi\|_{2,2} , 
$$
which gives (\ref{thm 2-0}). This completes the proof of Theorem~\ref{thm:max norm 1}.  
\end{proof} 

\subsection{The proof of Theorem~\ref{thm:convergence}}
Corollary~\ref{cor:max norm 2} is a direct consequence of Theorem~\ref{thm:max norm 1}, so that a uniform in time $\| \cdot \|_\infty$ bound of the numerical solution becomes available. With the help of the bound (\ref{cor 3-0}), we proceed into the convergence proof in Theorem~\ref{thm:convergence}.

\begin{proof} 
An application of the Taylor expansion for the exact solution $\Phi$ at $(x_i,y_j,z_k,t_{\kappa+\hf})$ implies that 
\begin{equation}
\begin{aligned}
&\frac{\Phi_{ijk}^{\kappa+1}-\Phi_{ijk}^\kappa}{\tau}=\trh U_{ijk}^{\kappa+\hf}+R_{ijk}^\kappa,\ \ 1 \leq i \leq m,\ \ 1 \leq j \leq n,\ \ 1 \leq k \leq l,\ \ 1 \leq s \leq N_t,\\
&U_{ij k}^{\kappa+\hf}=(\Phi_{ij k}^{\kappa+\hf})(\Phi_{ijk}^2)^{\kappa+\hf}+(1-\epsilon)\Phi_{ij k}^{\kappa+\hf}+3\trh\Phi_{ijk}^\kappa-\trh\Phi_{ijk}^{\kappa-1}+\trh(\trh\Phi_{ijk})^{\kappa+\hf}+S_{ijk}^\kappa,\\ & \ \ 1 \leq i \leq m,\ \ 1 \leq j \leq n,\ \ 1 \leq k \leq l,\ \ 1 \leq \kappa \leq N_t,\\
&\Phi_{ijk}^0=\psi(x_i,y_j,z_k)\ \ 1 \leq i \leq m,\ \ 1 \leq j \leq n,\ \ 1 \leq k \leq l , \label{1.30}
\end{aligned}
\end{equation}
in which the local truncation errors $R_{ijk}^\kappa $ and $S_{ijk}^\kappa$ satisfy 
\begin{equation} 
 |R_{ij k}^\kappa| \leq C_3(\tau^2+ h^2) , \, \, \, |S_{ij k}^\kappa| \leq C_3(\tau^2+ h^2) , \quad 1 \leq i \leq m,\ \ 1 \leq j \leq n,\ \ 1 \leq k \leq l,\ \ 1 \leq \kappa \leq N_t,\label{1.40} 
\end{equation}
for some $C_3\geq 0$, dependent only on $T$ and $L$. 

To facilitate the convergence analysis, we denote 
\begin{equation} 
  C_4= \| \Phi \|_{L^\infty (0,T; \Omega)} . \label{max bound-exact} 
\end{equation} 
The uniform in time $\| \cdot \|_\infty$ bound for the numerical solution is given by $\tilde{C}_0$, as defined as (\ref{cor 3-0}). These two bounds will be useful in the nonlinear error estimate. 

Subtracting \eqref{1.10} from \eqref{1.30} leads to the following error evolutionary equation: 
\begin{equation}
\begin{aligned}
\frac{e^{\kappa+1}-e^\kappa}{\tau}= \trh\left\{    {\left[(\Phi^{\kappa+\hf})(\Phi^2)^{\kappa+\hf}-(\phi^{\kappa+\hf})(\phi^2)^{\kappa+\hf}\right]+(1-\epsilon) e^{\kappa+\hf}}\right.\\
\left.{+(3\trh e^\kappa-\trh e^{\kappa-1})+\trh(\trh e)^{\kappa+\hf}+S^\kappa}   \right\}+R^\kappa .\label{1.50}
\end{aligned}
\end{equation}
Taking a discrete inner product with \eqref{1.50} by $e^{\kappa+\hf}$, we get 
\begin{equation}
\begin{aligned}
\left(\frac{e^{\kappa+1}-e^\kappa}{\tau} ,e^{\kappa+\hf}\right)&=\left(\Phi^{\kappa+\hf}(\Phi^2)^{\kappa+\hf}-\phi^{\kappa+\hf}(\phi^2)^{\kappa+\hf},\trh e^{\kappa+\hf}\right)+(1-\epsilon)(\trh e^{\kappa+\hf},e^{\kappa+\hf})\\
&+\left(\trh(3\trh e^{\kappa}-\trh e^{\kappa-1}), e^{\kappa+\hf}\right)+\left( \trh^3 e^{\kappa+\hf},e^{\kappa+\hf}\right)\\
&+(\trh S^\kappa,e^{\kappa+\hf})+(R^\kappa,e^{\kappa+\hf}) . \label{1.60}
\end{aligned}
\end{equation}
For the left hand side of \eqref{1.60}, the following identity is valid: 
\begin{equation} 
  \left(\frac{e^{\kappa+1}-e^\kappa}{\tau} ,e^{\kappa+\hf}\right)=\frac{1}{2\tau}(\|e^{\kappa+1}\|_2^2 - \|e^\kappa\|_2^2) . \label{error-1} 
\end{equation} 
For the nonlinear error term on the right hand side of \eqref{1.60}, we have 
\begin{equation} \label{error-2} 
\begin{aligned}
&\quad \left(\Phi^{\kappa+\hf}(\Phi^2)^{\kappa+\hf}-\phi^{\kappa+\hf}(\phi^2)^{\kappa+\hf},\trh e^{\kappa+\hf}\right)\\
&=\left(\Phi^{\kappa+\hf}\frac{(\Phi^\kappa)^2+(\Phi^{\kappa+1})^2}{2}-\phi^{\kappa+\hf}\frac{(\phi^\kappa)^2+(\phi^{\kappa+1})^2}{2},\trh e^{\kappa+\hf}\right)\\
&=\left(\Phi^{\kappa+\hf}\left(\frac{(\Phi^{\kappa+1})^2-(\phi^{\kappa+1})^2}{2}+\frac{(\Phi^\kappa)^2-(\phi^\kappa)^2}{2}\right)+(\Phi^{\kappa+\hf}-\phi^{\kappa+\hf})\frac{(\phi^\kappa)^2+(\phi^{\kappa+1})^2}{2},\trh e^{\kappa+\hf}\right)\\
&=\left(\Phi^{\kappa+\hf}\left(\frac{\Phi^{\kappa+1}+\phi^{\kappa+1}}{2} e^{\kappa+1}+\frac{\Phi^\kappa+\phi^\kappa}{2} e^\kappa\right)+\frac{(\phi^\kappa)^2+(\phi^{\kappa+1})^2}{2} e^{\kappa+\hf},\trh e^{\kappa+\hf}\right)\\
&\leq \frac14 \left[2C_4^2+(C_4 + \tilde{C}_0)^2\right] ( \| e^{\kappa+1} \|_2 + \| e^\kappa \|_2 ) \cdot  \| \trh e^{\kappa+\hf} \|_2 \\
& \le \frac12 \| \trh e^{\kappa+\hf} \|_2^2 + C  ( C_4^4 + \tilde{C}_0^4 )  ( \| e^{\kappa+1} \|_2^2 + \| e^\kappa \|_2^2 ) , 
\end{aligned}
\end{equation}
in which the discrete H\"older inequality has been repeatedly applied. 
The second and fourth terms on the right hand side of \eqref{1.60} could be analyzed in a straightforward way: 
\begin{eqnarray} 
  &&
 (1-\epsilon)(\trh e^{\kappa+\hf},e^{\kappa+\hf})=-(1-\epsilon)\|\nh e^{\kappa+\hf}\|_2^2 ,   \label{error-3} 
\\
  &&
  \left( \trh^3 e^{\kappa+\hf},e^{\kappa+\hf}\right)=-\|\nh \trh e^{\kappa+\hf}\|_2^2 . \label{error-4}   
\end{eqnarray}
For the third term of the right hand side of \eqref{1.60}, the following estimate is applied: 
\begin{eqnarray} \label{error-5} 
&&\quad \left(\trh(3\trh e^{\kappa}-\trh e^{\kappa-1}), e^{\kappa+\hf}\right)\\
&=&(3\trh e^{\kappa}-\trh e^{\kappa-1}, \trh e^{\kappa+\hf}) =\left(3\trh e^{\kappa}-\trh (2e^{\kappa-\hf}-e^\kappa), \trh e^{\kappa+\hf}\right)\\
&=&(4\trh e^{\kappa}-2\trh e^{\kappa-\hf}, \trh e^{\kappa+\hf})\\
&=&\left((2\trh e^{\kappa+1}+2\trh e^{\kappa})-(2\trh e^{\kappa+1}-2\trh e^{\kappa})-2\trh e^{\kappa-\hf},\trh e^{\kappa+\hf}\right)\\
&=& 4\|\trh e^{\kappa+\hf}\|_2^2 - \|\trh e^{\kappa+1}\|_2^2 + \|\trh e^\kappa\|_2^2 - 2(\trh e^{\kappa-\hf}, \trh e^{\kappa+\hf})\\
& \leq& 4\|\trh e^{\kappa+\hf}\|_2^2 - \|\trh e^{\kappa+1}\|_2^2 + \|\trh e^\kappa\|_2^2 + \| \trh e^{\kappa-\hf}\|_2^2 + \|\trh e^{\kappa+\hf}\|_2^2\\
& \leq& 5\|\trh e^{\kappa+\hf}\|_2^2 + \|\trh e^{\kappa-\hf}\|_2^2 - \|\trh e^{\kappa+1}\|_2^2 + \|\trh e^\kappa\|_2^2.
\end{eqnarray}
The fifth term and sixth terms on the right hand side of \eqref{1.60} could be bounded with an application of Cauchy inequality: 
\begin{eqnarray} 
  &&
  (\trh S^\kappa,e^{\kappa+\hf})\leq \hf \|S^\kappa\|_2^2+\hf \|\trh e^{\kappa+\hf}\|_2^2 , \label{error-6}  
\\
  &&
  (R^\kappa,e^{\kappa+\hf}) \leq \hf \|R^\kappa\|_2^2 + \hf \|e^{\kappa+\hf}\|_2^2 . \label{error-7}
\end{eqnarray}  
  
Therefore, a substitution of \eqref{error-1}-\eqref{error-6} into \eqref{1.60} yields 
\begin{equation} \label{error-8}
\begin{aligned}
&\quad \frac{1}{2\tau}(\|e^{\kappa+1}\|_2^2 - \|e^\kappa\|_2^2)+\|\nh \trh e^{\kappa+\hf}\|_2^2+(1-\epsilon)\|\nh e^{\kappa+\hf}\|_2^2\\
& \leq  C  ( C_4^4 + \tilde{C}_0^4 )  ( \| e^{\kappa+1} \|_2^2 + \| e^\kappa \|_2^2 ) +\hf \|e^{\kappa+\hf}\|_2^2 + \hf(\|S^\kappa\|_2^2 + \|R^\kappa\|_2^2) \\
&+6\|\trh e^{\kappa+\hf}\|_2^2 + \|\trh e^{\kappa-\hf}\|_2^2 - \|\trh e^{\kappa+1}\|_2^2 + \|\trh e^\kappa\|_2^2\\
& \le C  ( C_4^4 + \tilde{C}_0^4 +1 )  ( \| e^{\kappa+1} \|_2^2 + \| e^\kappa \|_2^2 ) +\hf(\|S^\kappa\|_2^2 + \|R^\kappa\|_2^2) \\
&+6\|\trh e^{\kappa+\hf}\|_2^2 + \|\trh e^{\kappa-\hf}\|_2^2 - \|\trh e^{\kappa+1}\|_2^2 + \|\trh e^\kappa\|_2^2 .
\end{aligned}
\end{equation}
By Lemma \ref{lemma2}, we obtain
\begin{eqnarray} 
  && 
  6 \|\trh e^{\kappa+\hf}\|_2^2 \leq 6\left(\frac{12^2}{3}\|e^{\kappa+\hf}\|_2^2 + \frac{2}{36}\|\nh \trh e^{\kappa+\hf}\|_2^2 \right) , 
  \label{error-9-1} 
\\
  &&
  \|\trh e^{\kappa-\hf}\|_2^2 \leq \left(\frac{1^2}{3}\|e^{\kappa-\hf}\|_2^2 + \frac{2}{3}\|\nh \trh e^{\kappa-\hf}\|_2^2\right) . 
  \label{error-9-2} 
\end{eqnarray} 
Going back \eqref{error-8}, we arrive at 
\begin{equation}  \label{error-10}
\begin{aligned}
&\quad \frac{1}{2\tau}(\|e^{\kappa+1}\|_2^2 - \|e^\kappa\|_2^2)+\|\nh \trh e^{\kappa+\hf}\|_2^2 + \|\trh e^{\kappa+1}\|_2^2 - \|\trh e^\kappa\|_2^2  \\
& \leq  288 \|e^{\kappa+\hf}\|_2^2 + \frac{1}{3} \|\nh \trh e^{\kappa+\hf}\|_2^2 + \frac{1}{3}\|e^{\kappa-\hf}\|_2^2 + \frac{2}{3}\|\nh \trh e^{\kappa-\hf}\|_2^2 \\
& +\hf(\|S^\kappa\|_2^2 + \|R^\kappa\|_2^2) + C  ( C_4^4 + \tilde{C}_0^4 +1 )  ( \| e^{\kappa+1} \|_2^2 + \| e^\kappa \|_2^2 ) . \\
\end{aligned}
\end{equation}
A summation in time implies that 
\begin{eqnarray}  
  && 
  \frac{1}{2\tau}(\|e^{\kappa+1}\|_2^2 - \|e^0\|_2^2) + \|\trh e^{\kappa+1}\|_2^2 - \|\trh e^0\|_2^2   \nn 
\\
  &\le& 
  \hf \sum_{i=0}^\kappa (\|S^i\|_2^2 + \|R^i\|_2^2) + C  \sum_{i=0}^\kappa ( C_4^4 + \tilde{C}_0^4 +1 )  ( \| e^{i+1} \|_2^2 + \| e^i \|_2^2 + \| e^{i-1} \|_2^2 ) . \label{error-11}  
\end{eqnarray}
In turn, an application of discrete Gronwall inequality yields the desired convergence result (\ref{convergence-0}). This completes the proof of Theorem~\ref{thm:convergence}. 
\end{proof} 
\section{Nonlinear multigrid solvers}   \label{sec:multigrid}

In this section we present the details of the nonlinear multigrid method that we use for solving the semi-implicit numerical scheme \eqref{1.10}. The fully finite-difference scheme \eqref{1.10} is formulated as follows: Find $\phi_{i,j,k}^{\kappa+1}$, $\mu_{i,j,k}^{\kappa+\nicefrac{1}{2}}$ and $\kappa_{i,j,k}^{\kappa+\nicefrac{1}{2}}$ in ${\mathcal C}_{\rm per}$ such that

\begin{eqnarray} 
\phi_{i,j,k}^{\kappa+1}  - \tau d_x\left(M\left( A_x\phi^{\kappa+1/2}_*\right)D_x\mu^{\kappa+1/2}\right)_{i,j,k}&&\nonumber\\
-\tau d_y\left(M\left( A_y\phi^{\kappa+1/2}_*\right)D_y\mu^{\kappa+1/2}\right)_{i,j,k}
-\tau d_z\left(M\left(A_z\phi^{\kappa+1/2}_*\right)D_z\mu^{\kappa+1/2}\right)_{i,j,k} &=& \phi_{i,j,k}^\kappa  ,  \label{2nd scheme-PFC-s1}
\\
\mu_{i,j,k}^{\kappa+1/2}  - \frac14 \left(\phi_{i,j,k}^{\kappa+1}+\phi_{i,j,k}^\kappa\right) \left((\phi_{i,j,k}^{\kappa+1})^2+(\phi_{i,j,k}^\kappa)^2\right)
- \frac{1-\epsilon}{2} \left(\phi_{i,j,k}^{\kappa+1}+\phi_{i,j,k}^\kappa\right)\nn&&\\
-3\Delta_h\phi_{i,j,k}^\kappa+\Delta_h\phi_{i,j,k}^{\kappa-1}-\Delta_h\omega_{i,j,k}^{\kappa+1}&=& 0,\label{2nd scheme-PFC-s2}
\\
\omega_{i,j,k}^{\kappa+1/2}-\frac{1}{2}\left(\Delta_h\phi_{i,j,k}^{\kappa+1}+\Delta_h\phi_{i,j,k}^\kappa\right)&=& 0, 
\label{2nd scheme-PFC-s3}
\end{eqnarray}
where $\phi^{\kappa+1/2}_*=\frac32 \phi^\kappa - \frac12 \phi^{\kappa-1}$. Denote
${\bf u}=(\phi_{i,j,k}^{\kappa+1}, \mu_{i,j,k}^{\kappa+\nicefrac{1}{2}},\omega_{i,j,k}^{\kappa+\nicefrac{1}{2}})^T$. Then the above discrete nonlinear system can be written in terms of a nonlinear operator ${\bf N}$ and source term ${\bf S}$ such that 
\begin{eqnarray}\label{eqn:dis-nonlinear-system}
{\bf N}({\bf u})= {\bf S} . 
\end{eqnarray} 
The $3\times m\times n\times l$ nonlinear operator ${\bf N}({\bf u}^{\kappa+1})=\big(N_{i,j,k}^{(1)}({\bf u}), N_{i,j,k}^{(2)}({\bf u}), N_{i,j,k}^{(3)}({\bf u})\big)^T$ can be defined as
\begin{eqnarray} 
N_{i,j,k}^{(1)}({\bf u})&=&\phi_{i,j,k}^{\kappa+1}  - \tau d_x\left(M\left( A_x\phi^{\kappa+1/2}_*\right)D_x\mu^{\kappa+1/2}\right)_{i,j,k}-
\tau d_y\left(M\left( A_y\phi^{\kappa+1/2}_*\right)D_y\mu^{\kappa+1/2}\right)_{i,j,k}\nn\\
&&-\tau d_z\left(M\left(A_z\phi^{\kappa+1/2}_*\right)D_z\mu^{\kappa+1/2}\right)_{i,j,k},  \label{scheme-PFC-o1}
\\
N_{i,j,k}^{(2)}({\bf u})&=&\mu_{i,j,k}^{\kappa+1/2}  - \frac14 \left(\phi_{i,j,k}^{\kappa+1}+\phi_{i,j,k}^\kappa\right) \left((\phi_{i,j,k}^{\kappa+1})^2+(\phi_{i,j,k}^\kappa)^2\right)
- \frac{1-\epsilon}{2} \phi_{i,j,k}^{\kappa+1}
-\Delta_h\omega_{i,j,k}^{\kappa+1},\label{scheme-PFC-o2}
\\
N_{i,j,k}^{(3)}({\bf u})&=&\omega_{i,j,k}^{\kappa+1/2}-\frac{1}{2} \Delta_h\phi_{i,j,k}^{\kappa+1} , 
\label{scheme-PFC-o3}
\end{eqnarray}
and the $3\times m\times n\times l$ source ${\bf S} =\big(S_{i,j,k}^{(1)}, S_{i,j,k}^{(2)}, S_{i,j,k}^{(3)}\big)^T$ is given by 
\begin{eqnarray} 
S_{i,j,k}^{(1)}&=&\phi_{i,j,k}^\kappa  ,  \label{scheme-PFC-s1}
\\
S_{i,j,k}^{(2)} &=&\frac{1-\epsilon}{2} \phi_{i,j,k}^\kappa
+3\Delta_h\phi_{i,j,k}^\kappa-\Delta_h\phi_{i,j,k}^{\kappa-1},\label{scheme-PFC-s2}
\\
S_{i,j,k}^{(3)}&=&\frac{1}{2}\Delta_h\phi_{i,j,k}^\kappa.
\label{scheme-PFC-s3}
\end{eqnarray}

The system \eqref{eqn:dis-nonlinear-system} can be efficiently solved  using a nonlinear Full Approximation Scheme (FAS) multigrid method (Algorithm~\ref{algorithm:FAS}: $\nu_1$ and $\nu_2$ are pre-smoothing and post-smoothing steps, $\ell, L$ are the current level and coarsest levels, and ${\bf I}_\ell^{\ell+1}, {\bf I}_{\ell+1}^\ell$ are coarsening and interpolating operators, respectively). More details  can be found in \cite{trottenberg2000multigrid}. Since we are using a standard FAS V-cycle approach, as reported in earlier works  \cite{baskaran13a,feng2016bsam,hu09,guo16,wise10}, we only provide the details of nonlinear smoothing scheme. For smoothing operator, we use a nonlinear Gauss-Seidel method with Red-Black ordering. 

\begin{algorithm}
	\caption{Nonlinear Multigrid Method (FAS)}
	\label{algorithm:FAS}
\begin{algorithmic}[1]
	\State Given $ {\bf u}^0$
	\Procedure{FAS}{${\bf N}^0, {\bf u}^0, {\bf S}^0, \nu_1, \nu_2, \ell=0$}
	\While{ $residual>tolerance$}
	\State Pre-smooth: ${\bf u}^\ell:$=smooth(${\bf N}^\ell, {\bf u}^\ell, {\bf S}^\ell, \nu_1$) \Comment{nonlinear Gauss-Seidel method }
	\State Residual : $ {\bf r}^\ell={\bf S}^\ell-{\bf N}^\ell {\bf u}^\ell$
	\State Coarsening: ${\bf r}^{\ell+1} = {\bf I}_\ell^{\ell+1}{\bf r}^\ell, {\bf u}^{\ell+1}={\bf I}_\ell^{\ell+1}{\bf u}^\ell$
	\If{$\ell=L$}
	\State Solve: ${\bf N}^{\ell+1} {\bf v}^{\ell+1}={\bf N}^{\ell+1} {\bf u}^{\ell+1}+{\bf r}^{\ell+1}$ \Comment{Cramer's rule }
	\State Error: ${\bf e}^{\ell+1}=  {\bf v}^{\ell+1}- {\bf u}^{\ell+1}$
	\Else
	\State Recursion: FAS({${\bf N}^{\ell+1}, {\bf u}^{\ell+1}, {\bf S}^{\ell+1}, \nu_1, \nu_2, \ell+1$})
	\EndIf
	\State Correction: ${\bf u}^\ell ={\bf u}^\ell+ {\bf I}_{\ell+1}^\ell{\bf e}^{\ell+1}$
	\State Post-smooth:  ${\bf u}^\ell:$=smooth(${\bf N}^\ell, {\bf u}^\ell, {\bf S}^\ell, \nu_2$) \Comment{nonlinear Gauss-Seidel method}
	\EndWhile\EndProcedure
\end{algorithmic}
\end{algorithm}

Let $n$ be the smoothing iteration, and define
\begin{equation*}
\begin{aligned}
M_{i+\hf,j,k}^{x} &:= M\left(\frac{1}{2}A_x \phi_{i+\hf,j,k}^\kappa-\frac{1}{2}A_x \phi_{i+\hf,j,k}^{\kappa-1}\right),\\
M_{i,j+\hf,k}^{y} &:= M\left(\frac{1}{2}A_y \phi_{i,j+\hf,k}^\kappa-\frac{1}{2}A_y \phi_{i,j+\hf,k}^{\kappa-1}\right),\\
M_{i,j,k+\hf}^{z} &:= M\left(\frac{1}{2}A_z \phi_{i,j,k+\hf}^\kappa-\frac{1}{2}A_z \phi_{i,j,k+\hf}^{\kappa-1}\right).
\end{aligned}
\end{equation*}  
Then the smoothing scheme is given by: for every $(i,j,k)$, stepping  lexicographically form $(1,1,1)$ to $(m,n,l)$, find $\phi_{i,j,k}^{\kappa+1,n+1}, \mu_{i,j,k}^{\kappa+\hf,n+1}, \kappa_{i,j,k}^{\kappa+\hf,n+1}$ that solve 
\begin{equation*} 
\begin{aligned}
\phi_{i,j,k}^{\kappa+1,n+1} +\frac{\tau}{h^2}\left( M_{i+\hf,j,k}^{x} +M_{i-\hf,j,k}^{x}+M_{i,j+\hf,k}^{y}+M_{i,j-\hf,k}^{y}+M_{i,j,k+\hf}^{z}+M_{i,j,k-\hf}^{z}\right)\mu_{i,j,k}^{\kappa+\hf,n+1}&=\tilde{S}_{i,j,k}^{(1)},
\\
\mu_{i,j,k}^{\kappa+1/2,n+1}  - \left(\frac14\left((\phi_{i,j,k}^{\kappa+1,n})^2+(\phi_{i,j,k}^\kappa)^2\right)+\frac{1-\epsilon}{2}\right)\phi_{i,j,k}^{k+1,n+1} +\frac{6}{h^2}\omega_{i,j,k}^{\kappa+\hf,n+1}&= \tilde{S}_{i,j,k}^{(2)},
\\
\omega_{i,j,k}^{\kappa+1/2}+\frac{3}{h^2}\phi_{i,j,k}^{\kappa+1,n+1}&= \tilde{S}_{i,j,k}^{(3)}, 
\end{aligned}
\end{equation*}
where
\begin{equation*} 
\begin{aligned}
\tilde{S}_{i,j,k}^{(1)}:=& S_{i,j,k}^{(1)}+\frac{\tau}{h^2}\left( M_{i+\hf,j,k}^{x}\mu_{i+1,j,k}^{\kappa+\hf,n}+M_{i-\hf,j,k}^{x}\mu_{i-1,j,k}^{\kappa+\hf,n+1}+M_{i,j+\hf,k}^{y}\mu_{i,j+1,k}^{\kappa+\hf,n}+M_{i,j-\hf,k}^{y}\mu_{i,j-1,k}^{\kappa+\hf,n+1}\right.\nn\\
&\left.+M_{i,j,k+\hf}^{z}\mu_{i,j,k+1}^{\kappa+\hf,n}+M_{i,j,k-\hf}^{z}\mu_{i,j,k-1}^{\kappa+\hf,n+1}\right), 
\\
\tilde{S}_{i,j,k}^{(2)}:=& S_{i,j,k}^{(2)}+\frac14\left((\phi_{i,j,k}^{\kappa+1,n})^2+(\phi_{i,j,k}^k)^2\right)\phi_{i,j,k}^k
\nn\\
&+\frac{1}{h^2}\left(\omega_{i+1,j,k}^{\kappa+\hf,n}+\omega_{i-1,j,k}^{\kappa+\hf,n+1}+\omega_{i,j+1,k}^{\kappa+\hf,n}+\omega_{i,j-1,k}^{\kappa+\hf,n+1}+\omega_{i,j,k+1}^{\kappa+\hf,n}+\omega_{i,j,k-1}^{\kappa+\hf,n+1}\right),
\\
\tilde{S}_{i,j,k}^{(3)}:=&
S_{i,j,k}^{(3)}+ \frac{1}{2h^2}\left(\phi_{i+1,j,k}^{\kappa+1,n}+\phi_{i-1,j,k}^{\kappa+1,n+1}+\phi_{i,j+1,k}^{\kappa+1,n}+\phi_{i,j-1,k}^{\kappa+1,n+1}+\phi_{i,j,k+1}^{\kappa+1,n}+\phi_{i,j,k-1}^{\kappa+1,n+1}\right). 
\end{aligned}
\end{equation*}
The above linearized system, which comes from  a local Picard linearization of the cubic term in the Gauss-Seidel scheme, can be solved by Cramer's Rule. 
\section{Numerical results}    \label{sec:numerical results}

In this section, we perform some numerical simulations for the three-dimensional scheme (\ref{1.10}), to verify the theoretical results. 

\subsection{Convergence and complexity of the multigrid solver}
In this subsection we demonstrate the accuracy and efficiency of the multigrid solver. We present the results of the convergence tests and perform some sample computations to verify the convergence and near optimal complexity with respect to the grid size $h$.  

In the first part of this test, we demonstrate the second order accuracy in time and space. The initial data is given by
\begin{eqnarray}\label{eqn:init1}
\phi_0(x,y,z) = 0.2+0.05\cos\big({2\pi x}/{3.2}\big)\cos\big({2\pi y}/{3.2}\big)\cos\big({2\pi z}/{3.2}\big),
\end{eqnarray} 
with $\Omega=[0, 3.2]^3$, $\epsilon = 2.5\times 10^{-2} $, $\tau=0.05h$ and $T=0.16$. We use a linear refinement path, \emph{i.e.}, $s=Ch$. At the final time $T=0.16$, we expect the global error to be $\mathcal{O}(s^2)+\mathcal{O}(h^2)=\mathcal{O}(h^2)$ under either the  $\ell^2$ or $\ell^\infty$ norm, as $h, s\to 0$.  Since we do not have an exact solution, instead of calculating the error at the final time, we compute the Cauchy difference, which is defined as $\delta_\phi: =\phi_{h_f}-\mathcal{I}_c^f(\phi_{h_c})$, where $\mathcal{I}_c^f$ is a bilinear interpolation operator (We applied Nearest Neighbor Interpolation in Matlab, which is similar to the 2D case in \cite{feng2016fch,feng2016preconditioned}). This requires having a relatively coarse solution, parametrized by $h_c$, and a relatively fine solution, parametrized by $h_f$, where $h_c = 2 h_f$, at the same final time. The $\ell^2$ norms of Cauchy difference and the convergence rates can be found in Table~\ref{tab:cov}. The  results confirm our expectation for the convergence order.

\begin{table}[!htb]
\begin{center}
\caption{Errors and convergence rates. Parameters are given in the
		text, and the initial data are defined in \eqref{eqn:init1}. The refinement path is $\tau=0.05h$.} \label{tab:cov}
\begin{tabular}{cccc}
\hline Grid sizes&$16^3- 32^3$ &$32^3- 64^3$& $64^3- 128^3$  \\
\hline Error& $ 2.3371\times 10^{-8}$&$ 5.8027 \times 10^{-9}$ &$ 1.4411 \times 10^{-9}$ 
\\Rate&-& 2.0099  &   2.0096   
\\
\hline
\end{tabular}
\end{center}
\end{table}

\begin{figure}[ht]
	\begin{center}
		\begin{subfigure}{0.45\textwidth}
			\includegraphics[width=\textwidth]{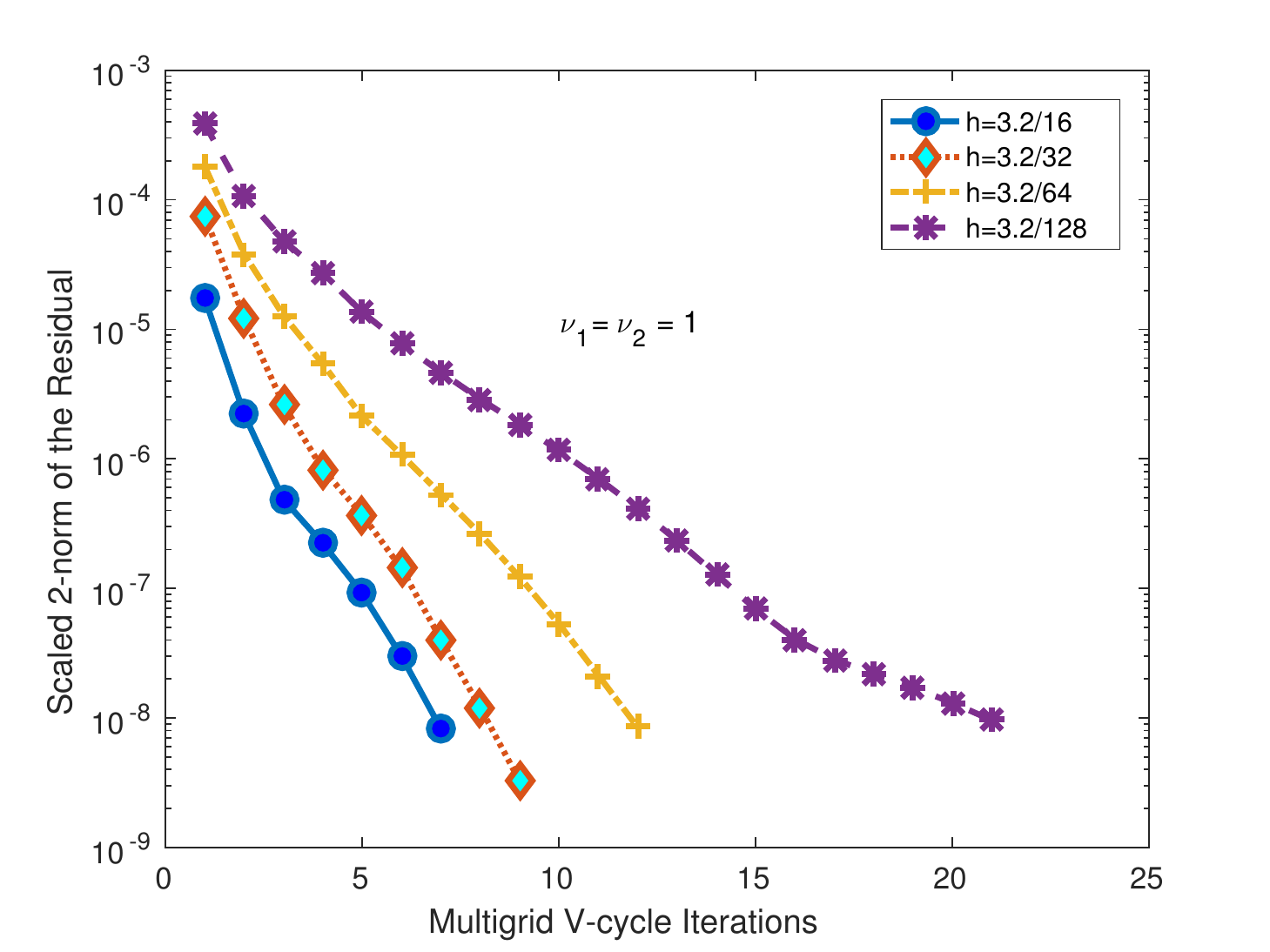} 
			\caption{$\nu_1=\nu_2 = 1$}
		\end{subfigure}
		\begin{subfigure}{0.45\textwidth}
			\includegraphics[width=\textwidth]{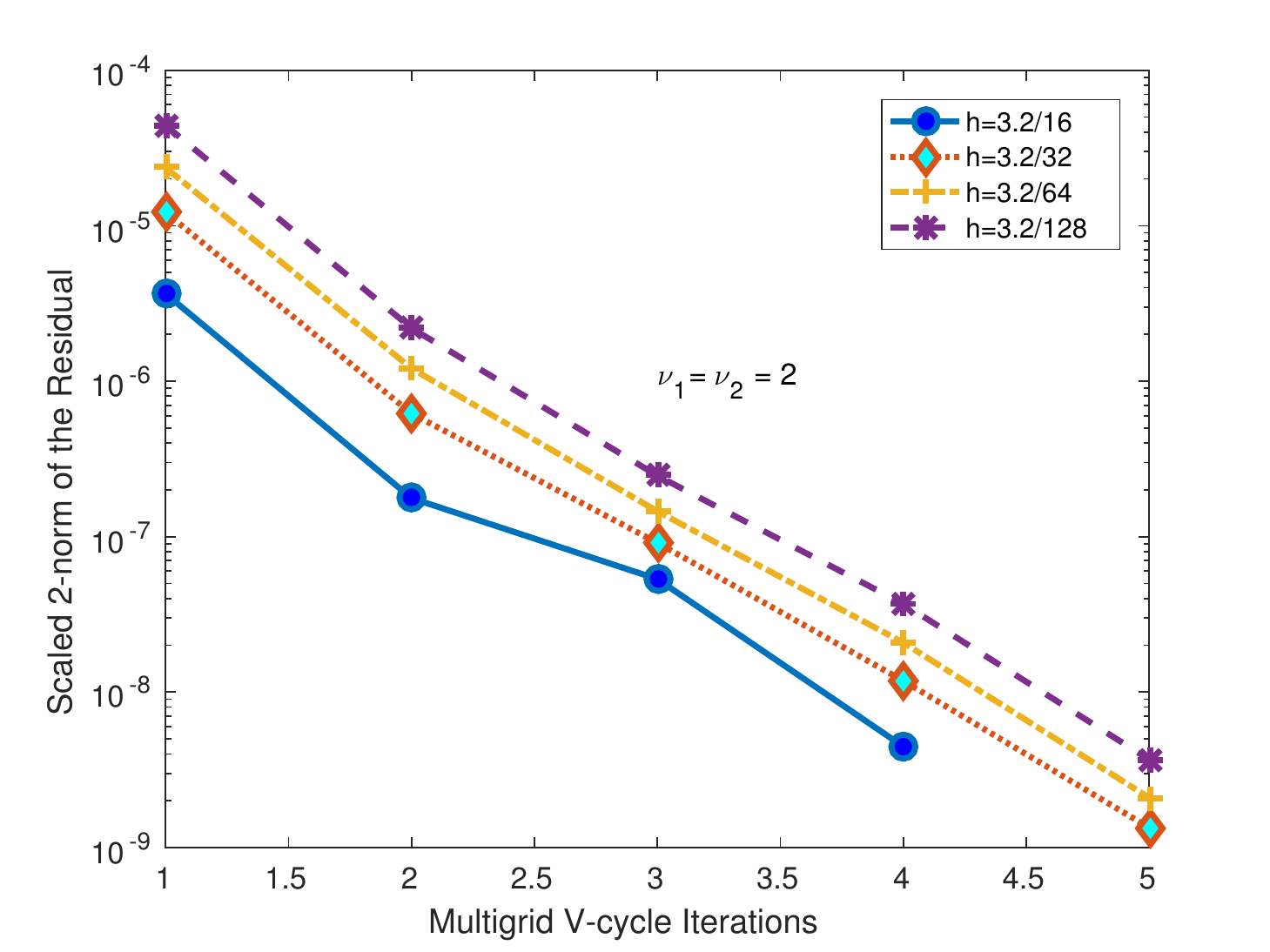}
			\caption{$\nu_1=\nu_2 = 2$}
		\end{subfigure}
	\end{center}
	\caption{The reduction in the norm of the residual for each V-cycle iteration at the 10th time step (i.e. at time $1.0\times 10^{-1}$ with time steps $\tau =1.0 \times 10^{-3}$). The rest of parameters are $\varepsilon = 5.0 \times 10^{-2}$, and $ \Omega=[0,3.2]\times[0,3.2]\times[0,3.2]$ and the initial condition is \eqref{eqn:init1}. }
	\label{fig:ch_complexity}
\end{figure}

\begin{rmk}
When calculating the Cauchy difference between the two different grids, the interpolation operator should be consist with the discrete stencil, otherwise the optimal convergence rate may not be observed. We applied Nearest Neighbor Interpolation in Matlab.
\end{rmk}
In the second part of this test, we investigate the complexity of the multigrid solver. The number of multigrid iterations to reach the residual tolerance is given in Table~\ref{tab:complexity}, for various choices of grid sizes and $(\nu_1, \nu_2)$. Table~\ref{tab:complexity}  indicates that the iteration numbers are nearly independent on $h$, when we use $2$ pre-smoothing and $2$ post-smoothing. The detailed reduction in the norm of the residual for each V-cycle iteration at the 10th time step can be found in Figure~\ref{fig:ch_complexity}. As can be
seen, the norm of the residual of each V-cycle is reduced by approximately the same rate each time with $\nu_1=\nu_2=2$, regardless of h. This is a typical feature of multigrid when it is operating with optimal complexity \cite{kay06, trottenberg2000multigrid, wise10}. For $\nu_1=\nu_2 = 1$, we do not observe a similar feature. Moreover, we also observe that more multigrid iterations are required for smaller values of h, which confirms our convergence analysis.
\begin{table}[!htb]
\begin{center}
\caption{The number of multigrid iterations of each residual below the tolerance $tol=10^{-8}$ at the 10-th time step (i.e. at time $1.0\times 10^{-2}$ with time steps $\tau =1.0 \times 10^{-3}$), The rest of parameters are $\varepsilon = 2.5 \times 10^{-2}$, and $ \Omega=[0,3.2]\times[0,3.2]\times[0,3.2]$.} \label{tab:complexity}
\begin{tabular}{ccccccccc}
\hline
&$(\nu_1,\nu_2)$\textbackslash Grid sizes& 16 & 32 & 64 & 128 \\
\hline &(1,1)& 7 & 9 & 12 & 21
\\ & $(2,2)$ & 4 & 5 & 5 & 5
\\
\hline
\end{tabular}
\end{center}
\end{table}

\subsection{Growth of a polycrystal}

\begin{figure}[!htp]
	\begin{center}
		\begin{subfigure}{0.49\textwidth}
			\includegraphics[width=0.491\textwidth,height=0.49\textwidth]{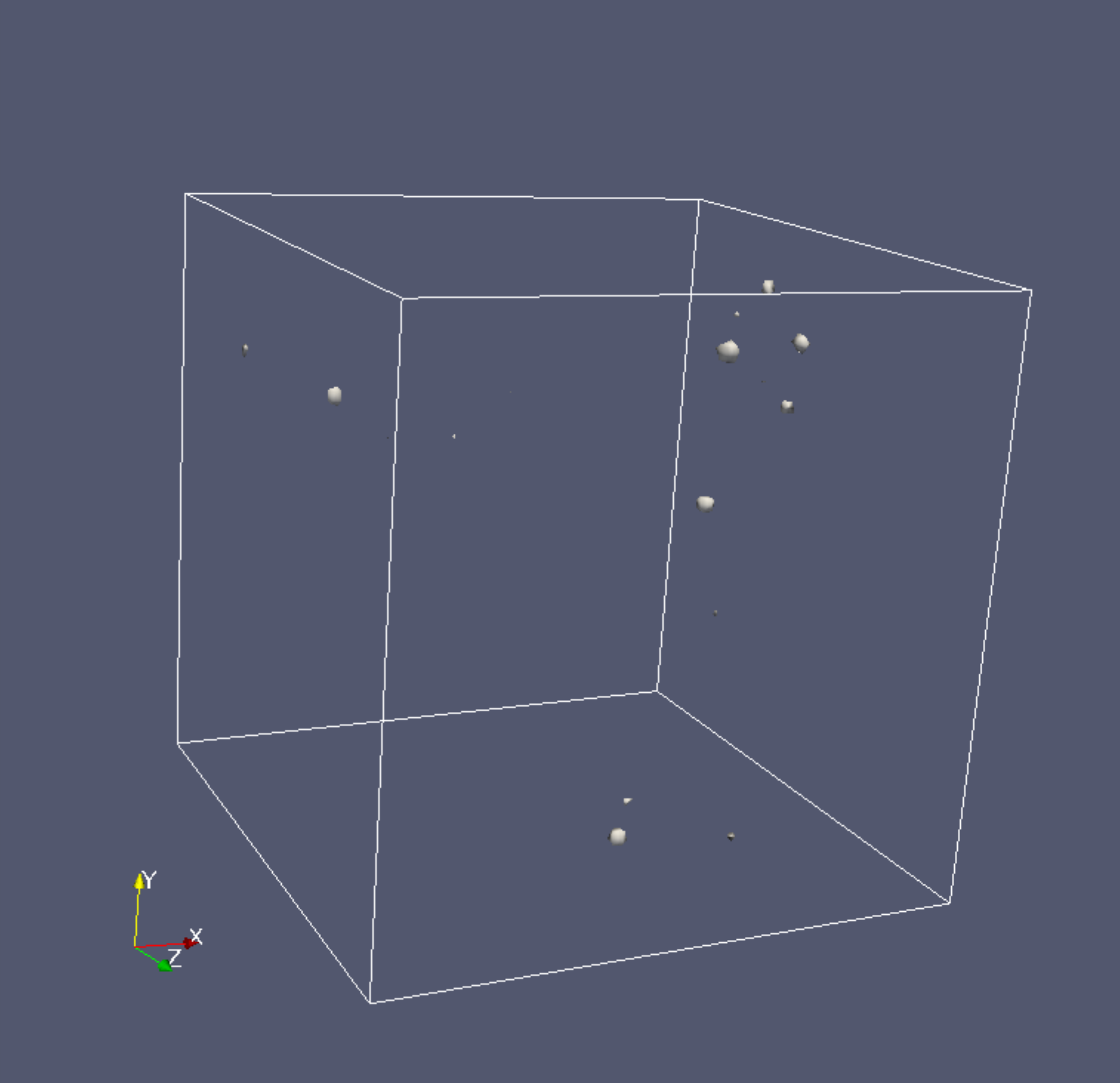} 
			\includegraphics[width=0.491\textwidth,height=0.49\textwidth]{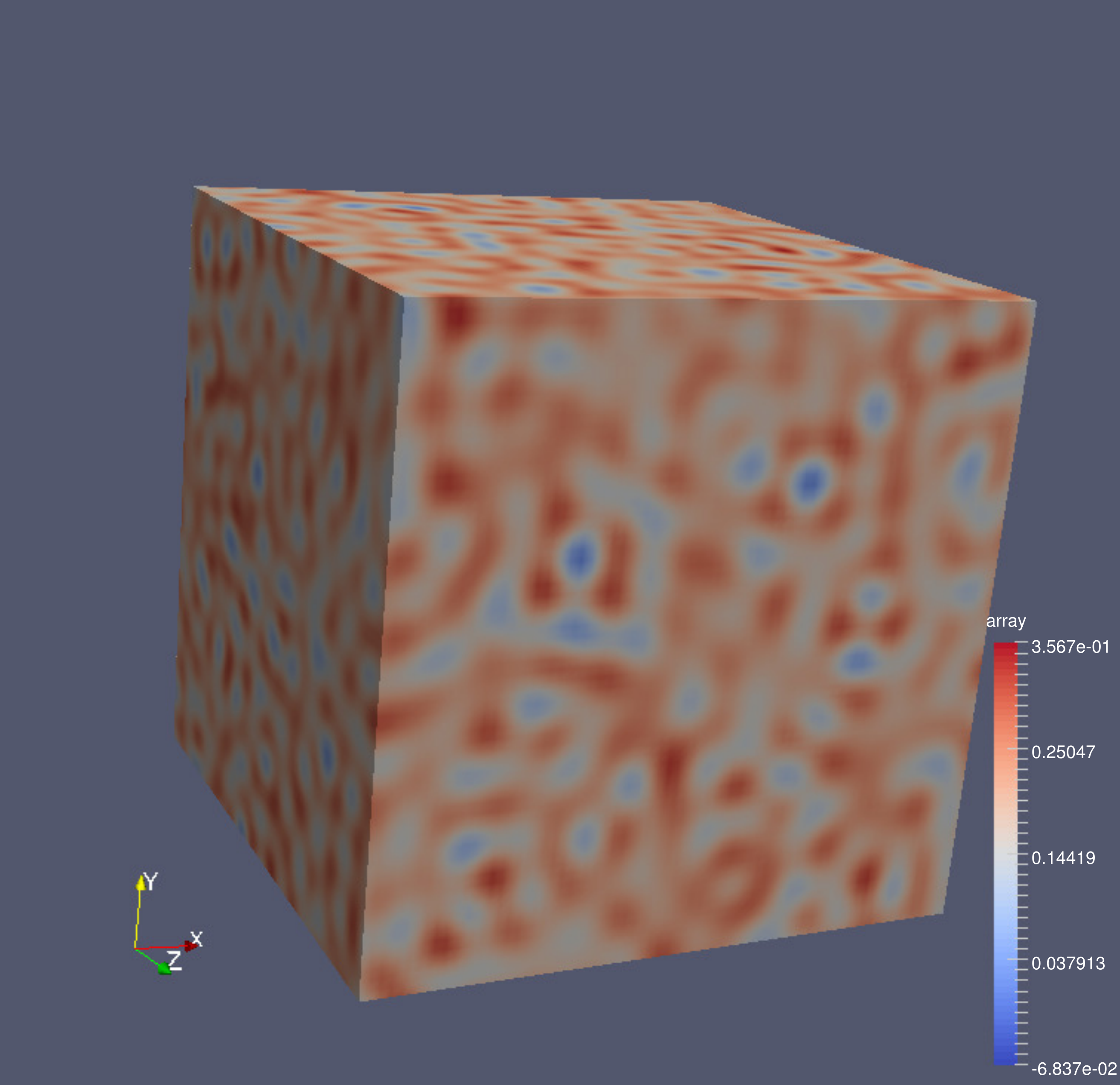}
			\caption*{$t=40$}
		\end{subfigure}
		\begin{subfigure}{0.49\textwidth}
			\includegraphics[width=0.49\textwidth,height=0.49\textwidth]{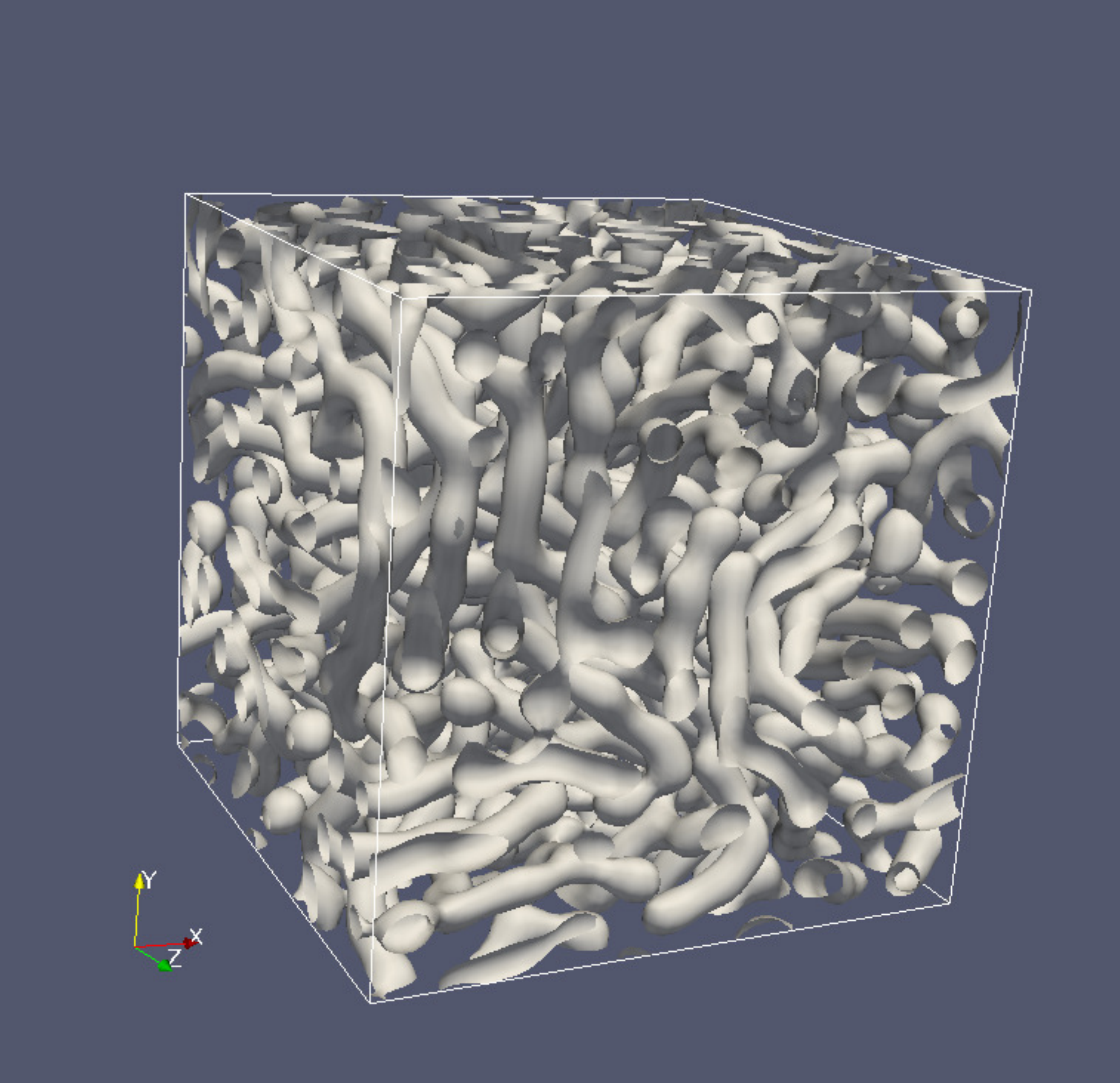} 
			\includegraphics[width=0.49\textwidth,height=0.49\textwidth]{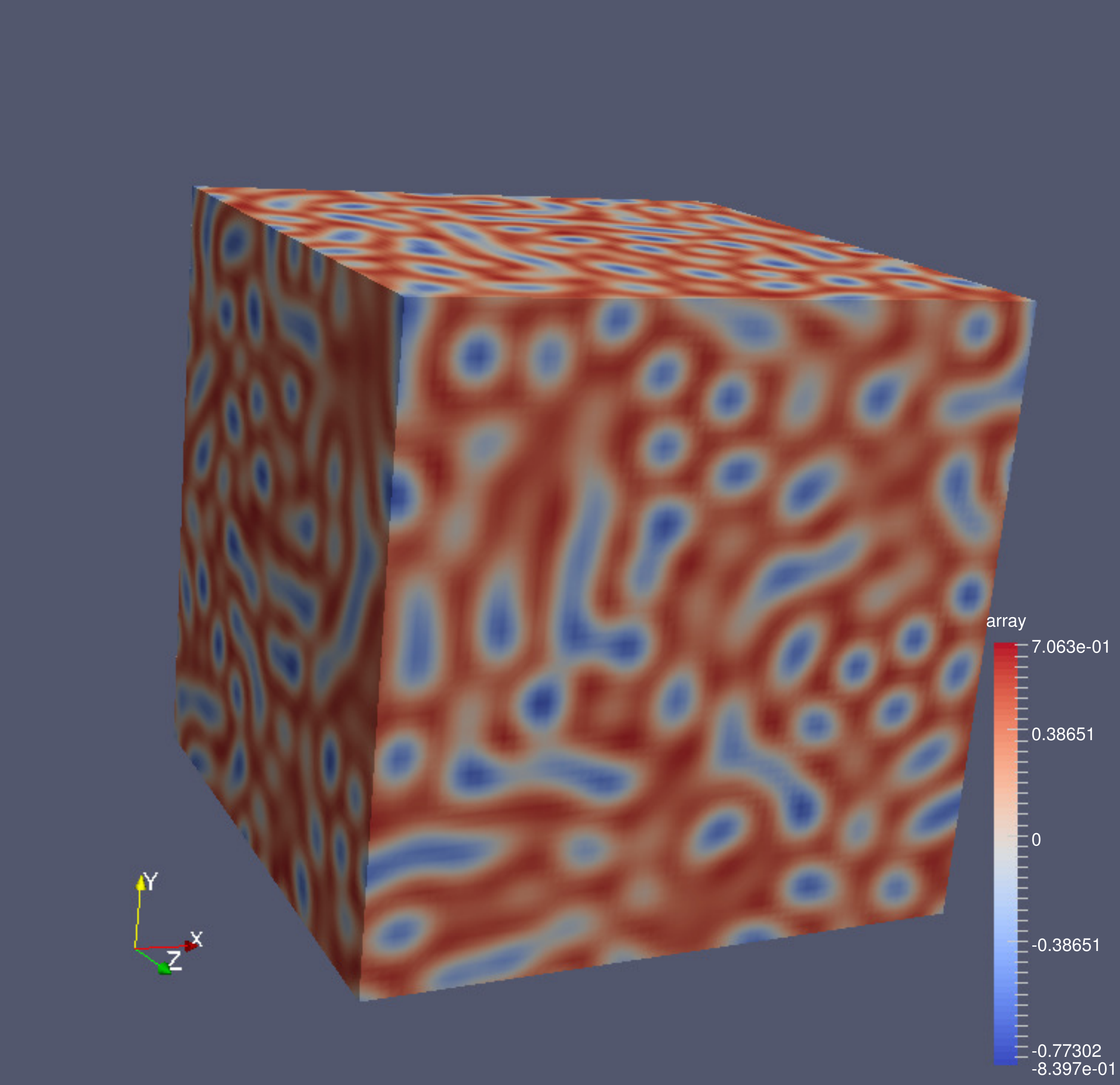}
			\caption*{$t=200$}
		\end{subfigure}
		\begin{subfigure}{0.49\textwidth}
			\includegraphics[width=0.49\textwidth,height=0.49\textwidth]{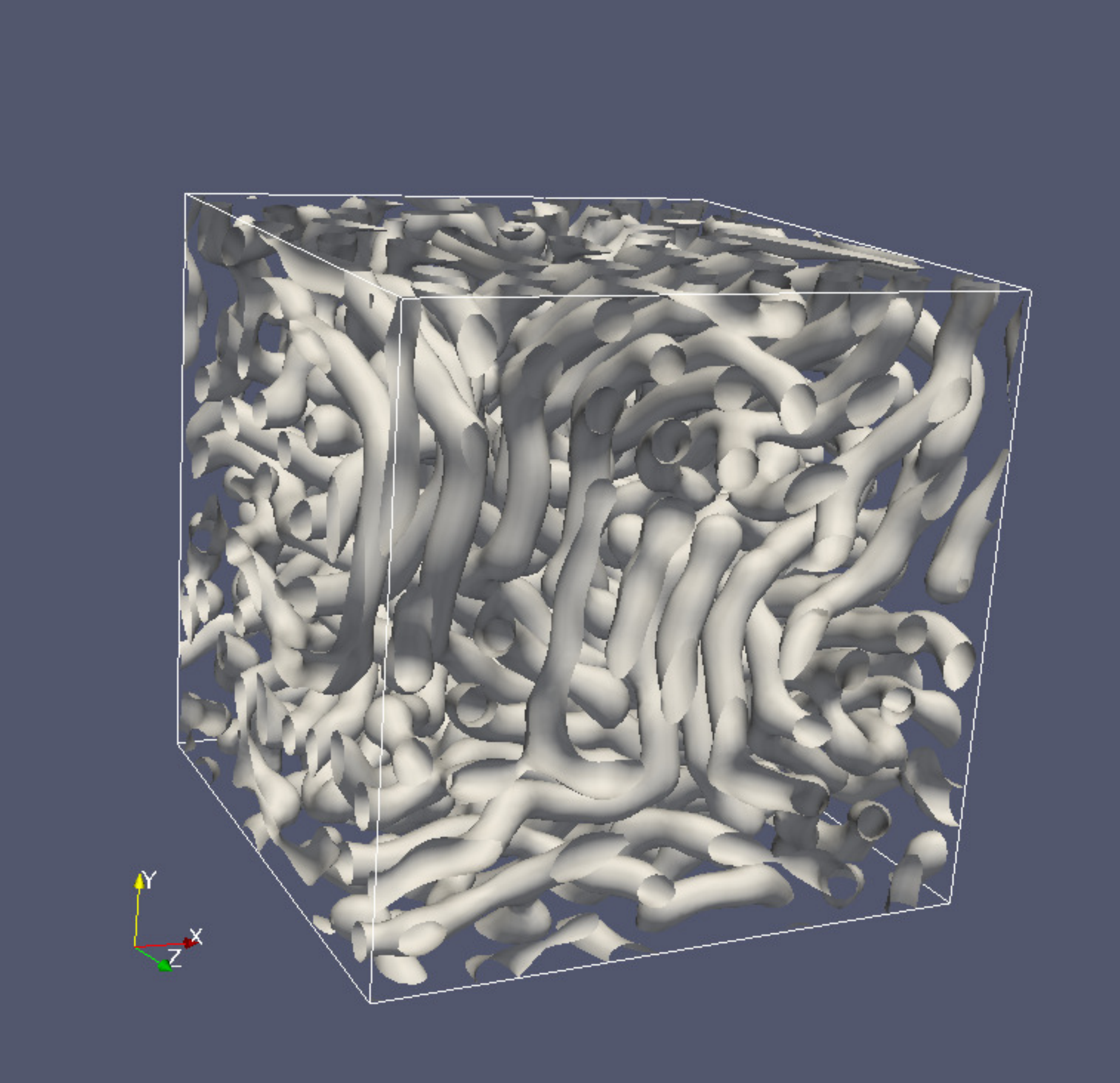} 
			\includegraphics[width=0.49\textwidth,height=0.49\textwidth]{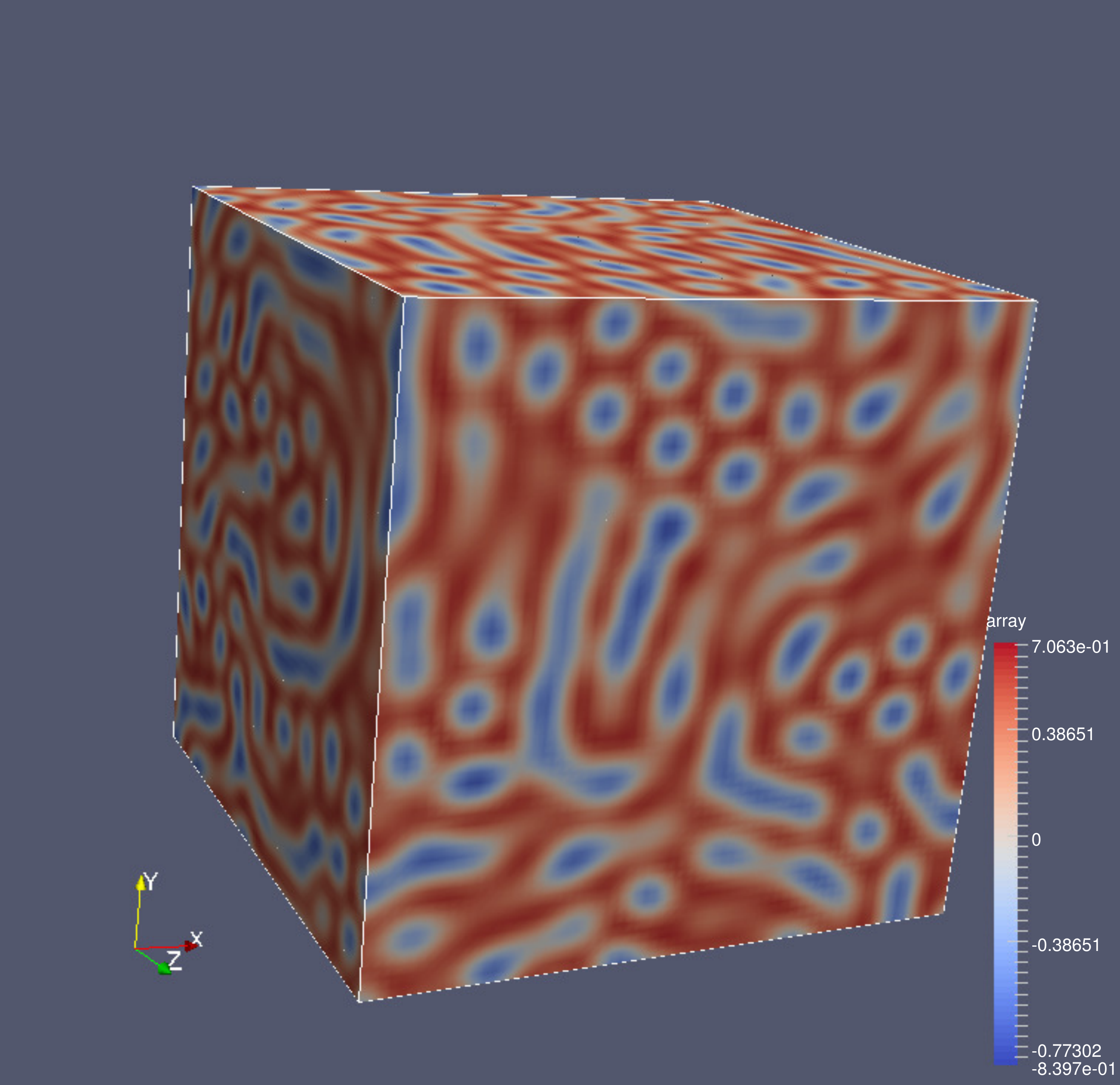}
			\caption*{$t=2000$}
		\end{subfigure}
		\begin{subfigure}{0.49\textwidth}
			\includegraphics[width=0.49\textwidth,height=0.49\textwidth]{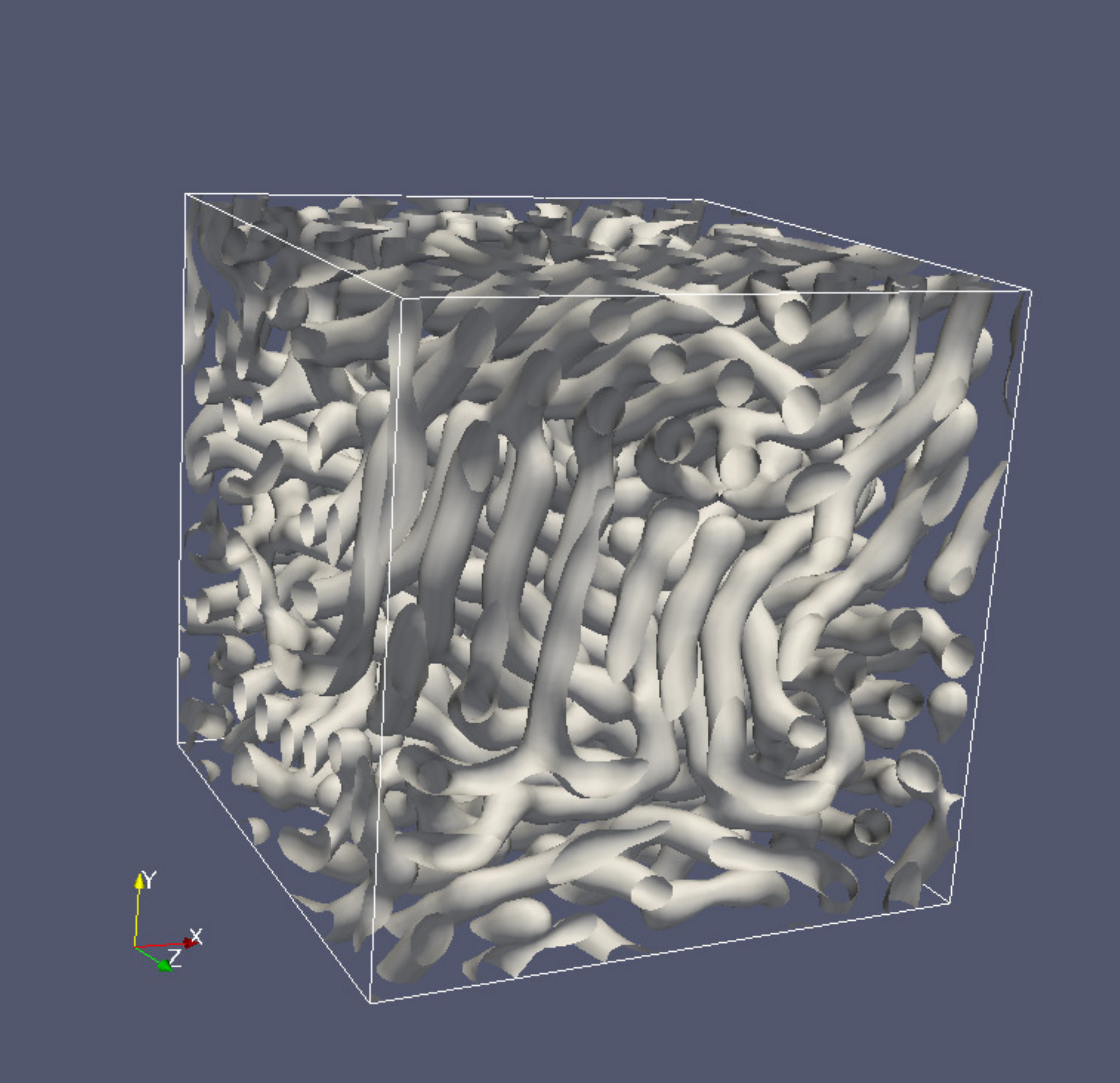} 
			\includegraphics[width=0.49\textwidth,height=0.49\textwidth]{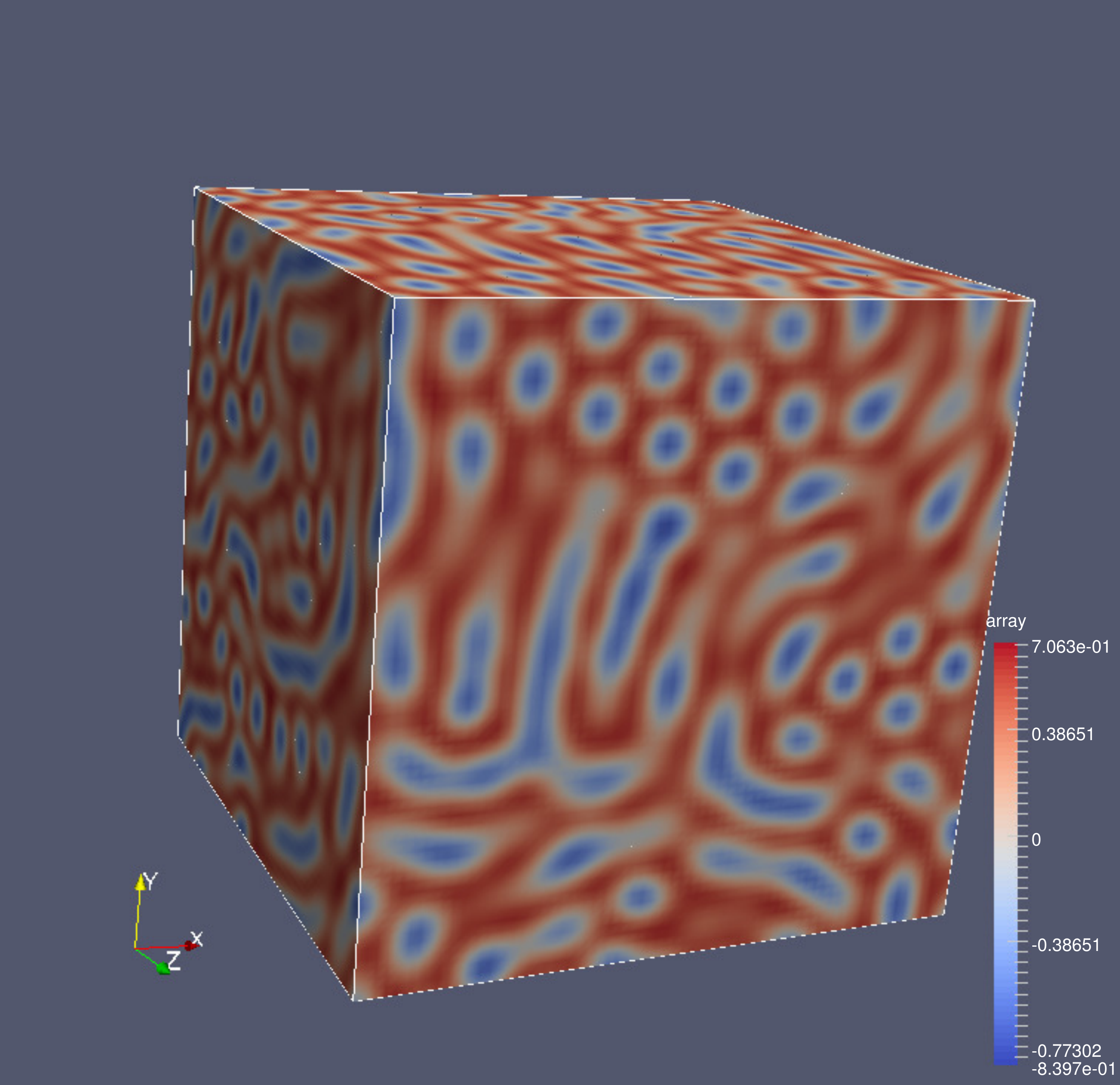}
			\caption*{$t=4000$}
		\end{subfigure}	
		\begin{subfigure}{0.49\textwidth}
			\includegraphics[width=0.49\textwidth,height=0.49\textwidth]{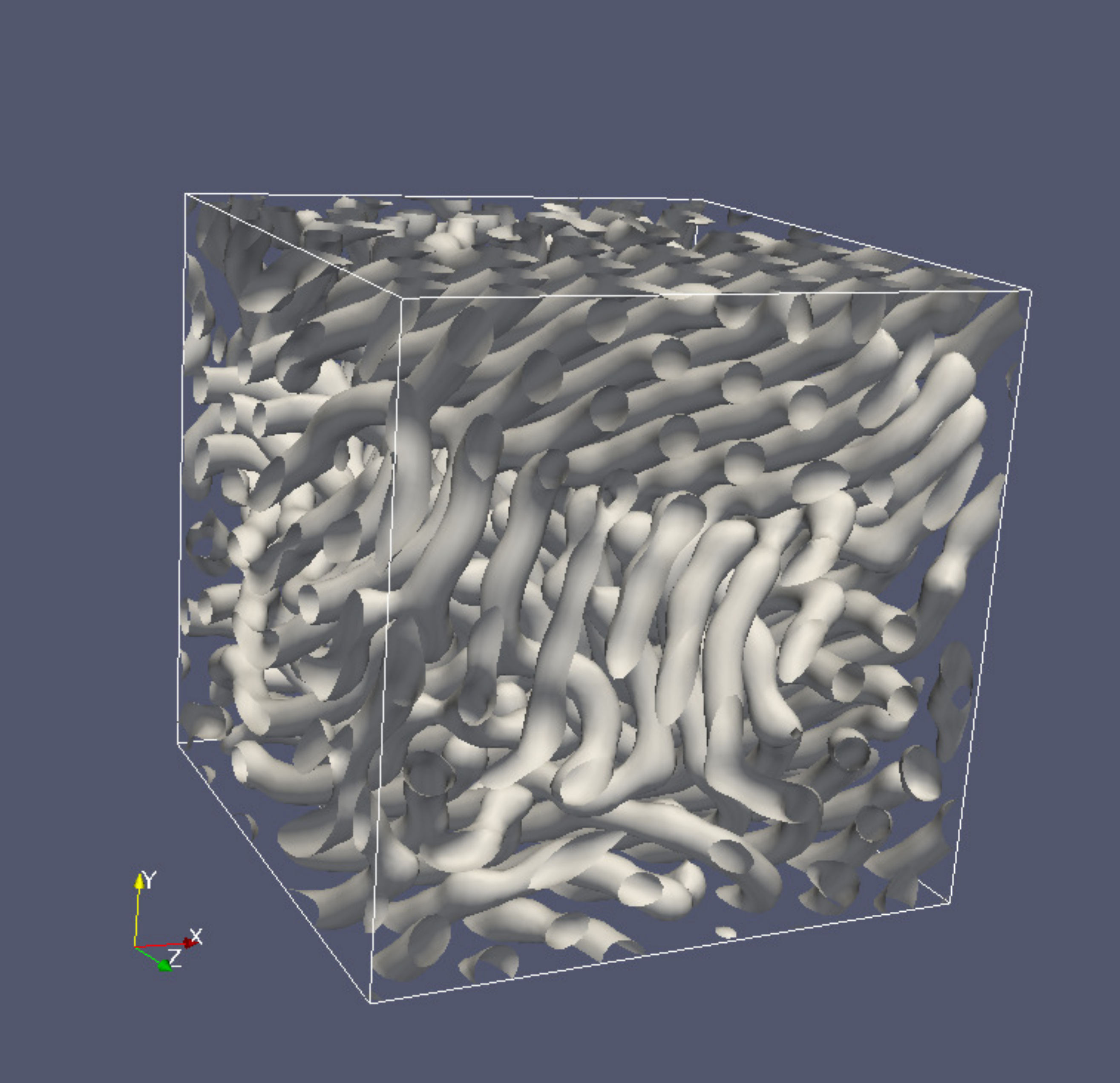} 
			\includegraphics[width=0.49\textwidth,height=0.49\textwidth]{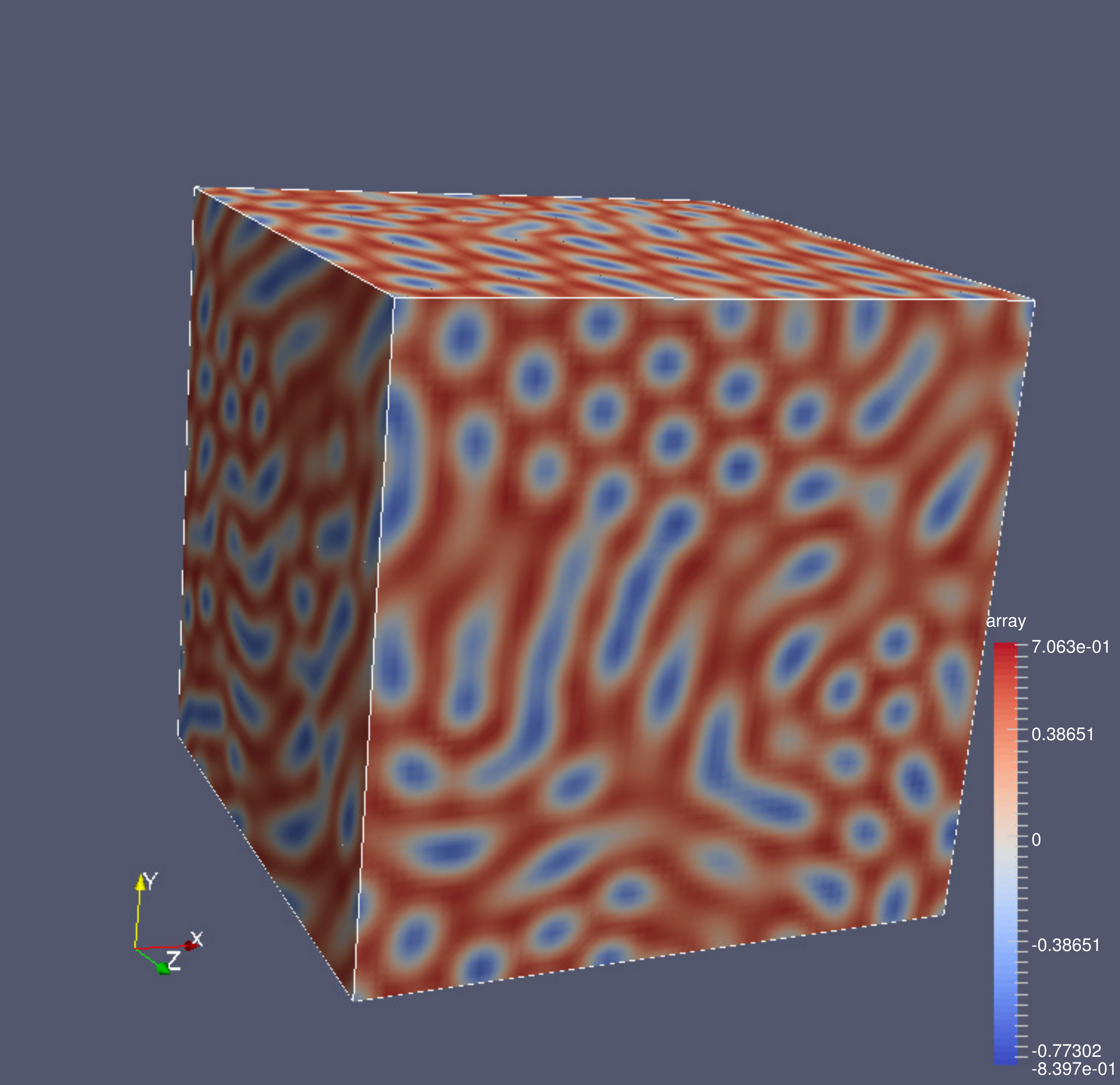}
			\caption*{$t=6000$}
		\end{subfigure}
		\begin{subfigure}{0.49\textwidth}
			\includegraphics[width=0.49\textwidth,height=0.49\textwidth]{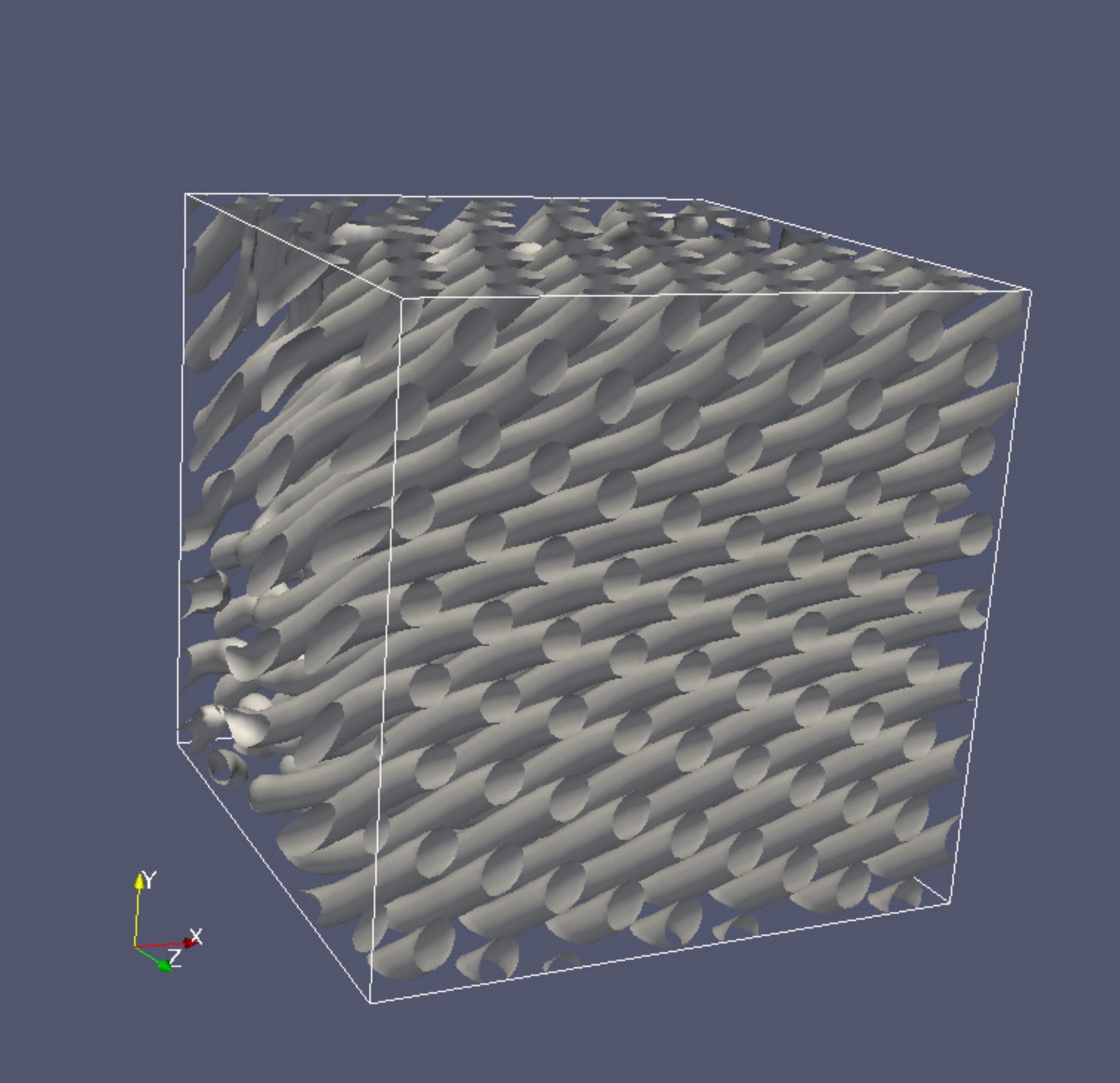} 
			\includegraphics[width=0.49\textwidth,height=0.49\textwidth]{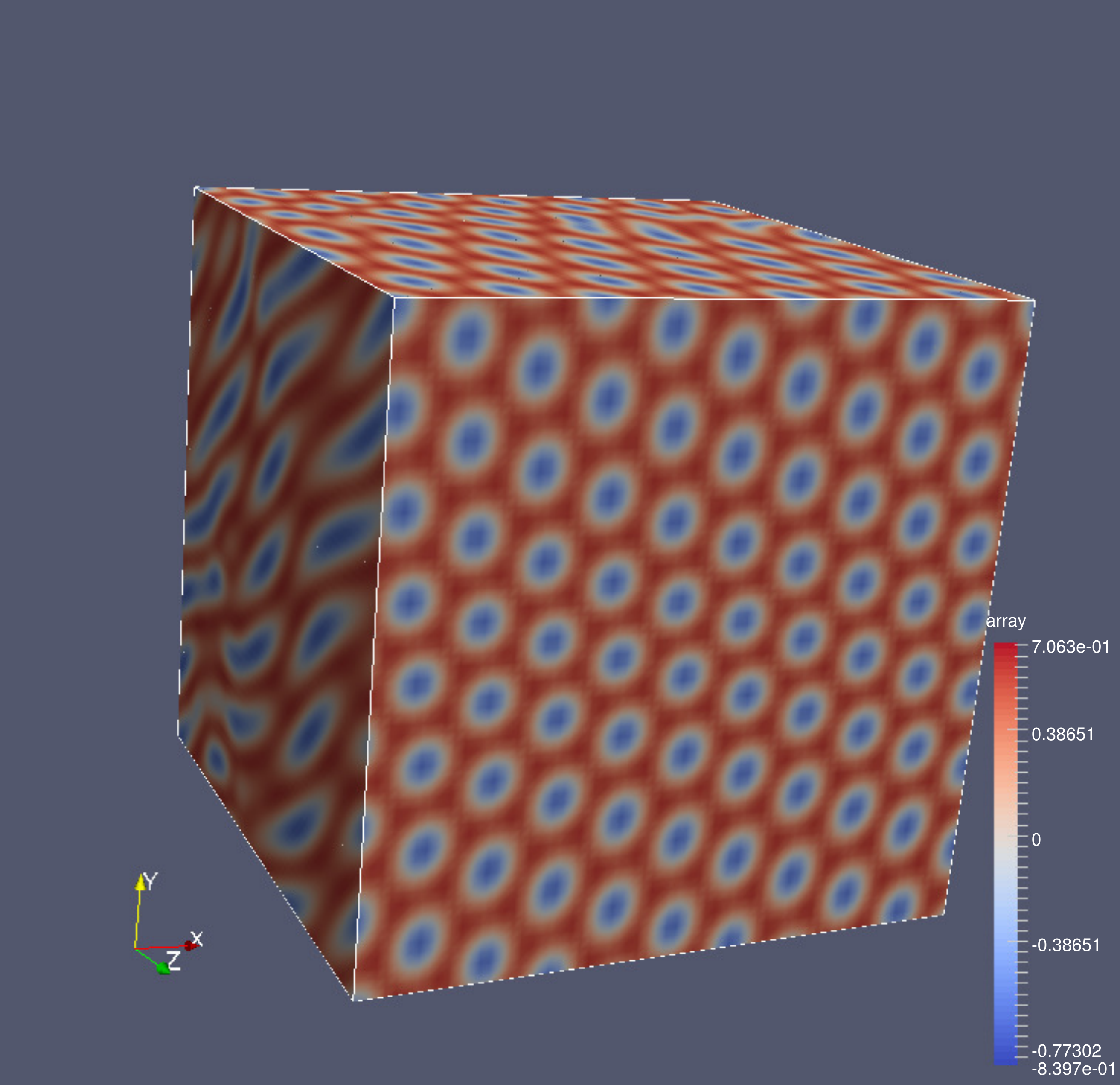}
			\caption*{$t=8000$}
		\end{subfigure}	
		\begin{subfigure}{0.49\textwidth}
			\includegraphics[width=0.49\textwidth,height=0.49\textwidth]{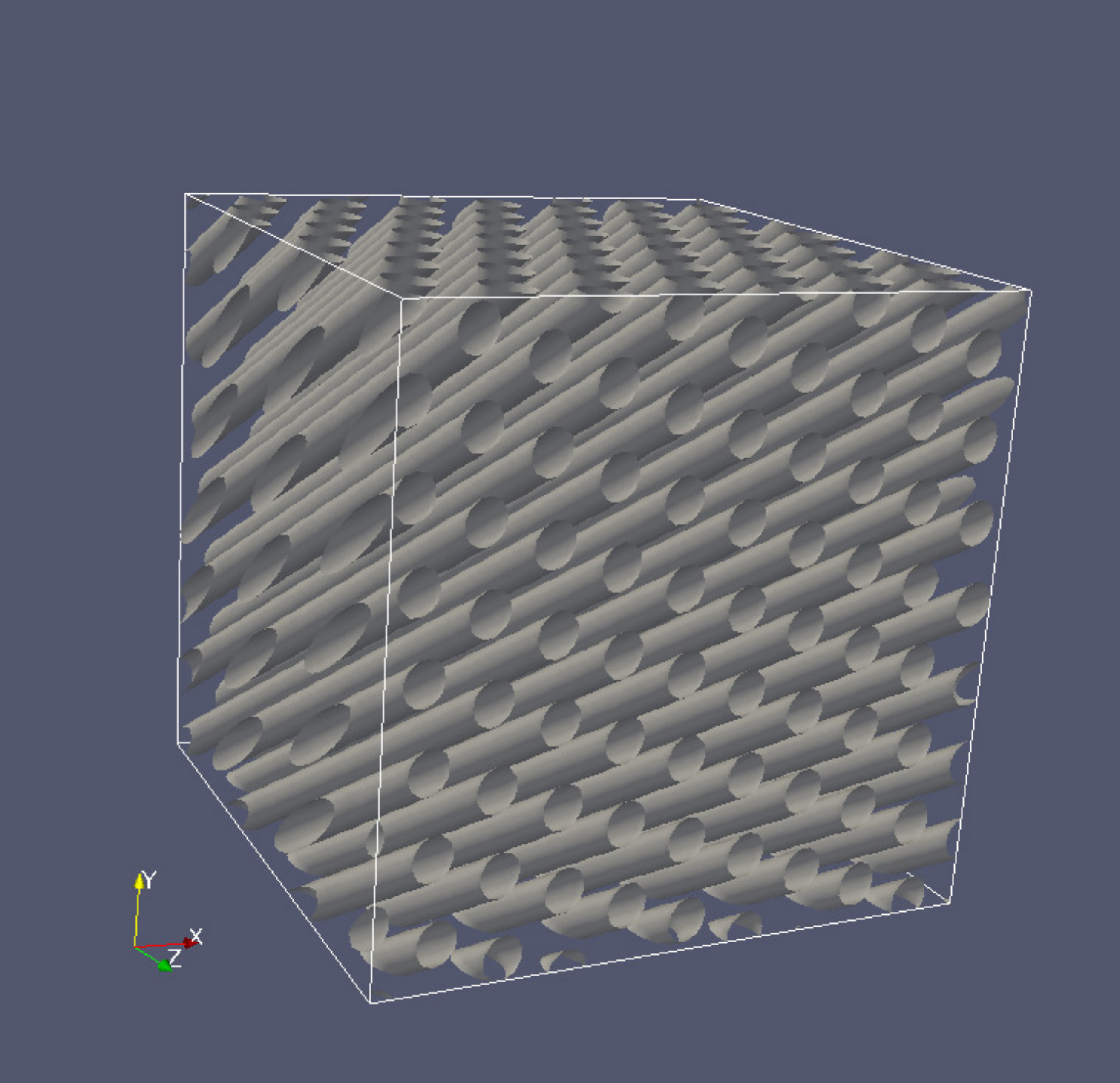} 
			\includegraphics[width=0.49\textwidth,height=0.49\textwidth]{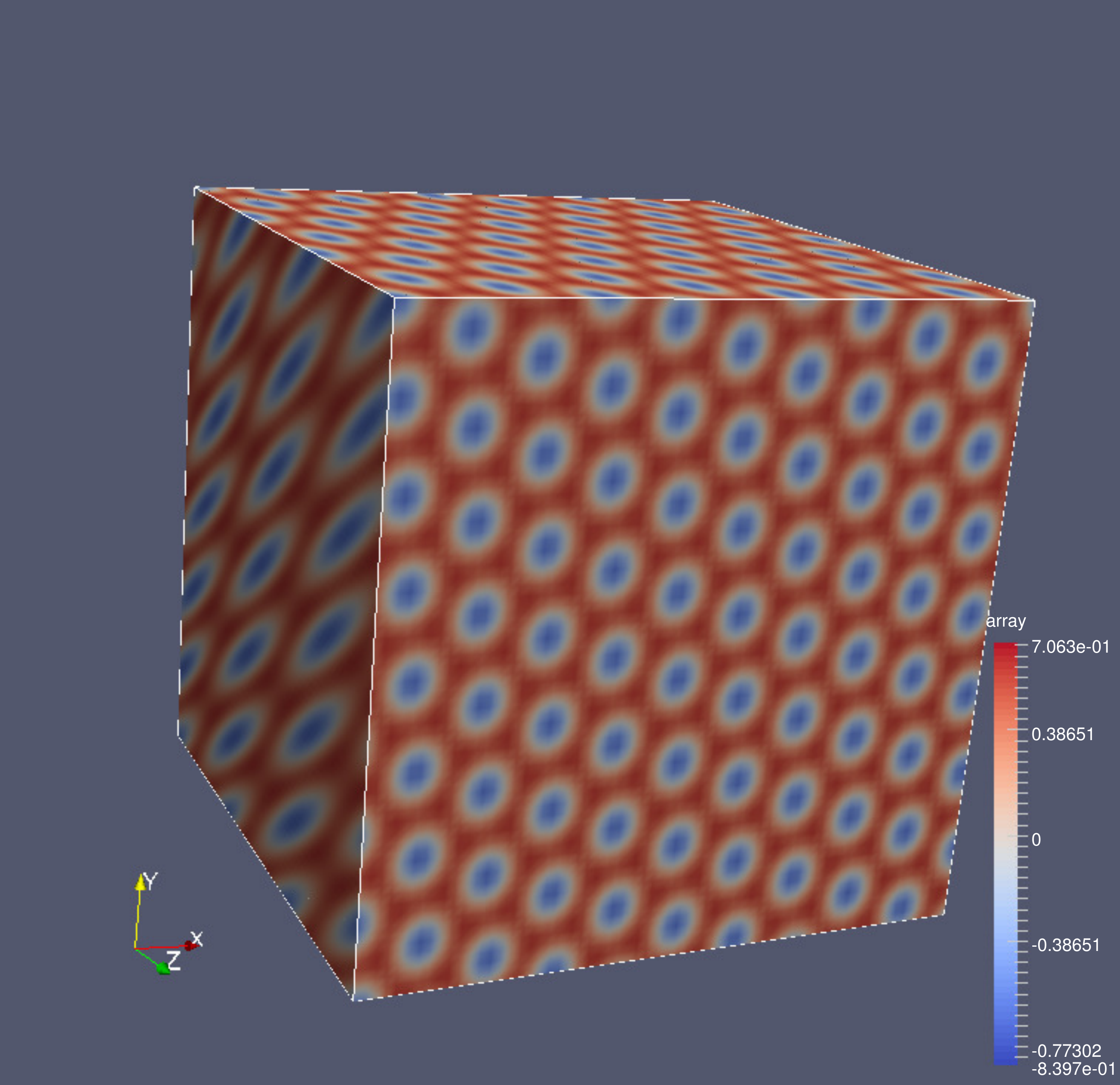}
			\caption*{$t=10000$}
		\end{subfigure}
		\begin{subfigure}{0.49\textwidth}
			\includegraphics[width=0.49\textwidth,height=0.49\textwidth]{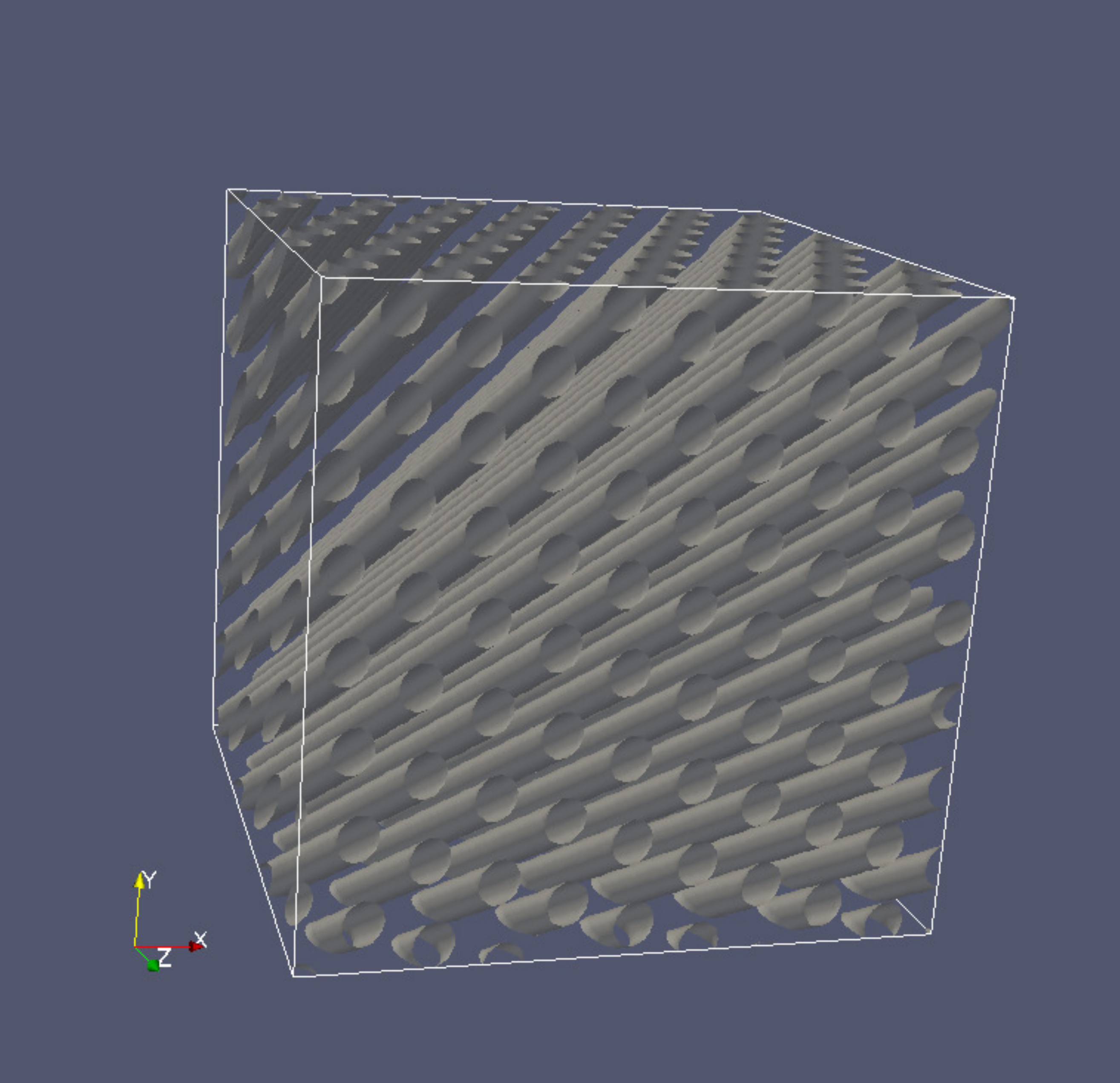} 
			\includegraphics[width=0.49\textwidth,height=0.49\textwidth]{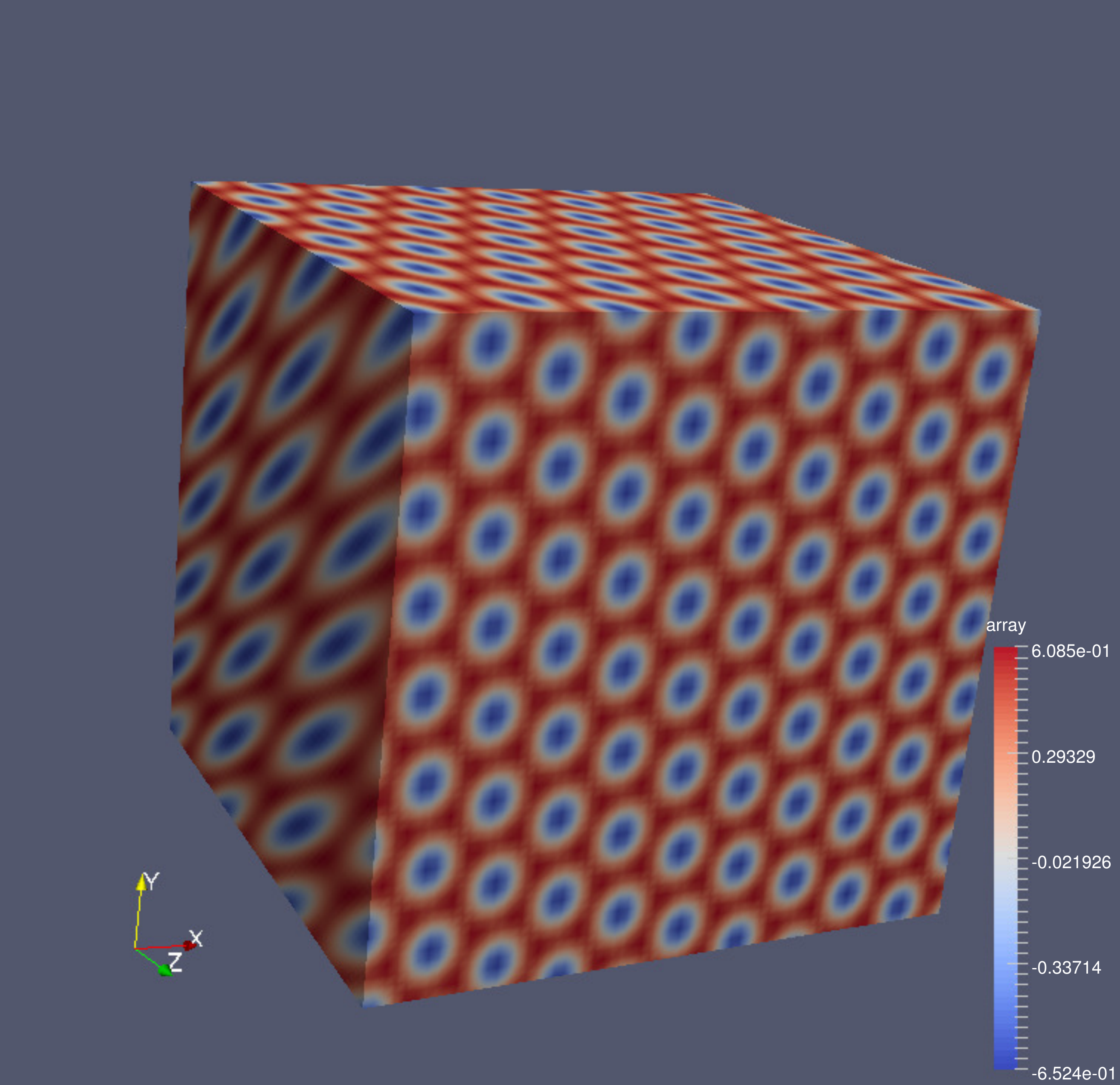}
			\caption*{$t=12000$}
		\end{subfigure}						
		\caption{Three-dimensional periodic micro-structures snapshots  with initial condition \eqref{eqn:init2}  at $t=40,200, 2000, 4000, 6000, 8000, 10000~ \text{and}~ 12000$. Left: iso-surface plots of $ \phi=0.0 $, Right: snapshots of micro-structures plots. The parameters are
			$\varepsilon = 2.5 \times 10^{-1}, \Omega=[0,64]\times[0,64]\times[0,64], \tau= 1.0 \times 10^{-2}$.}
		\label{fig:long-time-pfc}
	\end{center}
\end{figure}

The initial data for this simulation are taken as essentially random:
	\begin{equation}
	\label{eqn:init2}
    \phi^0_{i,j,k}=0.2+0.005\cdot r_{i,j,k},
	\end{equation}
where the $r_{i,j,k}$ are uniformly distributed random numbers in [0, 1]. Time snapshots of the micro-structure can be found in Figure \ref{fig:long-time-pfc}. The numerical results are consistent with the experiments on this topic in \cite{gomez2012unconditionally}.

\section{Conclusions}  \label{sec:conclusion} 
In this paper, we have provided a detailed convergence analysis of a finite difference scheme for the three-dimensional PFC equation, with the second order accuracy in both time and space established. The numerical scheme was proposed in \cite{hu09}, with the unique solvability and unconditional energy stability already proved in the earlier work. Meanwhile, a theoretical justification off the convergence analysis turns out to be challenging, due to a difficulty to obtain a maximum norm bound of the numerical solution in three-dimensional space. We overcome this difficulty with the help of discrete Fourier transformation, and repeated applications of Parseval equality in both continuous and discrete spaces. With such a discrete maximum norm bound developed for the numerical solution, the convergence analysis could be derived by a careful process of consistency estimate and stability analysis for the numerical error function. 

In addition, we describe the detailed multigrid solver to implement this numerical scheme over a three-dimensional domain. Various numerical results are presented, including the numerical convergence test and the three-dimensional polycrystal growth simulation. The efficiency and robustness of the nonlinear multigrid solver has been extensively demonstrated in these three-dimensional numerical experiments.

	\section{Acknowledgments}
The second author would like to thank Jing Guo at South China University of Technology for the valuable discussions. 
This work is supported in part by NSF DMS-1418689 (C.~Wang), NSF DMS-1418692 (S.~Wise),  NSFC 11271048, 91130021 and the Fundamental Research Funds for the Central Universities (Z.~Zhang).

\bibliographystyle{plain}
\bibliography{PFC3D}
\end{document}